\documentclass[preprint,3p,times]{elsarticle}

\usepackage{amsmath,amssymb,amsfonts}
\usepackage{bm}
\usepackage{graphicx}
\usepackage{booktabs}
\usepackage[section]{placeins}
\usepackage{siunitx}
\usepackage{subcaption}
\usepackage{tikz}
\usetikzlibrary{arrows.meta,positioning,calc,decorations.pathmorphing}
\usepackage{pgfplots}
\pgfplotsset{compat=1.16}
\pgfplotsset{colormap={viridis}{
  rgb=(0.267,0.005,0.329) rgb=(0.283,0.141,0.458) rgb=(0.254,0.265,0.530)
  rgb=(0.207,0.372,0.553) rgb=(0.164,0.471,0.558) rgb=(0.128,0.567,0.551)
  rgb=(0.135,0.659,0.518) rgb=(0.267,0.749,0.441) rgb=(0.478,0.821,0.318)
  rgb=(0.741,0.873,0.150) rgb=(0.993,0.906,0.144)}}
\pgfplotsset{colormap={magma}{
  rgb=(0.001,0.000,0.014) rgb=(0.078,0.054,0.211) rgb=(0.232,0.059,0.437)
  rgb=(0.390,0.070,0.503) rgb=(0.550,0.161,0.506) rgb=(0.716,0.215,0.475)
  rgb=(0.868,0.288,0.409) rgb=(0.967,0.440,0.360) rgb=(0.994,0.625,0.427)
  rgb=(0.996,0.812,0.573) rgb=(0.987,0.991,0.749)}}
\newcommand{\semicylcbar}[6]{%
\begin{minipage}[c]{0.15\linewidth}
\centering
\begin{tikzpicture}
\begin{axis}[
    hide axis, scale only axis, height=#5, width=1pt,
    colormap name=#1, colorbar, point meta min=0, point meta max=#2,
    colorbar style={
      width=0.26cm, height=#5,
      ytick={#3},
      yticklabel pos=right, ytick pos=right,
      ylabel={#4}, ylabel near ticks,
      tick label style={font=\footnotesize}, label style={font=\footnotesize},
      axis line style={line width=0.5pt, black}, tick align=outside,
    }]
\addplot [draw=none] coordinates {(0,0) (1,1)};
\end{axis}
\end{tikzpicture}
\end{minipage}}
\usepackage{listings}
\usepackage{xcolor}
\usepackage[hidelinks]{hyperref}

\definecolor{codegreen}{rgb}{0,0.5,0}
\definecolor{codegray}{rgb}{0.5,0.5,0.5}
\definecolor{codepurple}{rgb}{0.58,0,0.82}
\definecolor{backcolour}{rgb}{0.96,0.96,0.94}
\lstdefinestyle{python}{
  language=Python,
  backgroundcolor=\color{backcolour},
  commentstyle=\color{codegreen},
  keywordstyle=\color{blue},
  stringstyle=\color{codepurple},
  basicstyle=\ttfamily\footnotesize,
  breaklines=true, captionpos=b, keepspaces=true,
  showstringspaces=false, tabsize=2,
  frame=single, framesep=4pt, xleftmargin=6pt, xrightmargin=4pt,
}

\newcommand{\bphi}{\bm{\varphi}}
\newcommand{\bu}{\bm{u}}
\newcommand{\bA}{\bm{A}}
\newcommand{\bd}{\bm{d}}
\newcommand{\bn}{\bm{n}}
\newcommand{\bq}{\bm{q}}
\newcommand{\bw}{\bm{w}}
\newcommand{\bg}{\bm{g}}
\newcommand{\bdq}{\bm{dq}}
\newcommand{\bI}{\bm{\mathrm{id}}}
\newcommand{\bbeta}{\bm{\beta}}
\newcommand{\bT}{\bm{T}}
\newcommand{\Rbb}{\mathbb{R}}
\definecolor{shapeblue}{HTML}{0072B2}        

\journal{Computers \& Structures}

\newcommand{\areaBudgetBracket}{9}
\newcommand{\areaBudgetSemicyl}{2}

\newcommand{\bracketEvals}{119}
\newcommand{\bracketFallbacks}{3}

\newcommand{\bracketGainOverComsol}{5.6}
\newcommand{\bracketGrads}{84}
\newcommand{\bracketIters}{100}

\newcommand{\bracketReductionPct}{94.6}

\newcommand{\bracketWall}{1214}

\newcommand{\costAdjBracket}{1.20 [1.05, 1.26]}
\newcommand{\costAdjSemicyl}{0.30 [0.29, 1.07]}
\newcommand{\costContBracket}{43}
\newcommand{\costContSemicyl}{38}
\newcommand{\costFwdBracket}{9.64 [6.21, 39.14]}
\newcommand{\costFwdSemicyl}{1.36 [0.94, 2.45]}
\newcommand{\costGeomBracket}{0.019 [0.016, 0.022]}
\newcommand{\costGeomSemicyl}{0.010 [0.009, 0.011]}
\newcommand{\costRieszBracket}{0.009 [0.008, 0.013]}
\newcommand{\costRieszSemicyl}{0.004 [0.004, 0.005]}

\newcommand{\dispBudgetSemicyl}{0.1}

\newcommand{\firedrakeVersion}{2026.4.0}
\newcommand{\fireshapeCommitMain}{0340f76}
\newcommand{\fwdMaxErr}{0.0151}

\newcommand{\fwdRMSE}{0.0086}
\newcommand{\fwdRelErr}{0.50}
\newcommand{\hessianActions}{358}
\newcommand{\hessianGrads}{35}
\newcommand{\hessianQuality}{0.49}
\newcommand{\hessianReductionPct}{94.7}
\newcommand{\hessianReportedIter}{41}

\newcommand{\linFinalRawBracket}{70.762}

\newcommand{\linMaxDispBracket}{0.049}
\newcommand{\linMaxDispSemicyl}{0.094}
\newcommand{\linNonlinearRawBracket}{70.672}

\newcommand{\linOwnReductionPctBracket}{97.9}
\newcommand{\linOwnReductionPctSemicyl}{94.9}

\newcommand{\linRefFactorSemicyl}{16.3}
\newcommand{\linRefRatioBracket}{3.2}

\newcommand{\meshAngleBracket}{22}
\newcommand{\meshAngleSemicyl}{33}
\newcommand{\meshQualityBracket}{0.61}
\newcommand{\meshQualitySemicyl}{0.68}

\newcommand{\semicylAreaPct}{1.8}
\newcommand{\semicylAvgDeflPct}{87}
\newcommand{\semicylEvals}{116}
\newcommand{\semicylFallbacks}{17}

\newcommand{\semicylGrads}{87}
\newcommand{\semicylIters}{100}
\newcommand{\semicylMaxDisp}{0.095}
\newcommand{\semicylReductionPct}{98.4}

\newcommand{\semicylWall}{983}
\newcommand{\sinAmplitudeBracket}{5}
\newcommand{\sinAmplitudeListBracket}{$A = \pm 5, \pm 10, \pm 20$\,mm}
\newcommand{\sinAmplitudeListSemicyl}{$A = \pm 2, \pm 5, \pm 10, \pm 20, \pm 40$\,mm}
\newcommand{\sinAmplitudeSemicyl}{40}
\newcommand{\sinBudgetPctBracketSteps}{10, 20, 40}
\newcommand{\sinBudgetPctSemicylSteps}{2, 5, 10, 20, 40}

\newcommand{\sinModesLength}{6}
\newcommand{\sinModesWidth}{7}
\newcommand{\sinModesWordsBracket}{six half-waves along the L profile and seven across the width}
\newcommand{\sinModesWordsSemicyl}{six half-waves along the axis and seven across the arc}
\newcommand{\sinOppositeSignBracket}{-}
\newcommand{\sinOppositeSignSemicyl}{+}
\newcommand{\sinStartSignBracket}{+}
\newcommand{\sinStartSignSemicyl}{-}
\newcommand{\sinTicksBracket}{0,1.25,2.5,3.75,5}
\newcommand{\sinTicksSemicyl}{0,10,20,30,40}
\newcommand{\taylorRateBracket}{2.00}
\newcommand{\taylorRateSecondBracket}{2.96}
\newcommand{\taylorRateSecondSemicyl}{2.99}
\newcommand{\taylorRateSemicyl}{1.97}
\newcommand{\unboundedBracketAreaPct}{27.6}
\newcommand{\unboundedBracketFinalAreaPct}{44.8}

\newcommand{\unboundedBracketFinalReductionPct}{98.0}
\newcommand{\unboundedBracketReductionPct}{97.0}
\newcommand{\unboundedBracketReportedIter}{53}

\begin{document}

\begin{frontmatter}

\title{Shape optimisation of nonlinear Naghdi shells on
discrete geometries}

\author[1,2,3]{Ado Farsi\corref{cor1}}
\ead{ado.farsi@imperial.ac.uk}
\ead{ado.farsi@tanuki.ai}
\author[4]{Alberto Paganini}
\ead{a.paganini@leicester.ac.uk}
\cortext[cor1]{Corresponding author.}
\affiliation[1]{organization={Imperial College London, Department of Earth
Science and Engineering}, city={London}, country={United Kingdom}}
\affiliation[2]{organization={University College London, Department of Earth
Sciences}, city={London}, country={United Kingdom}}
\affiliation[3]{organization={Tanuki Technologies}, city={London},
country={United Kingdom}}
\affiliation[4]{organization={University of Leicester, School of Computing and Mathematical Sciences},
city={Leicester}, country={United Kingdom}}

\begin{abstract}
We present an automated framework for the shape optimisation of a geometrically nonlinear,
five-parameter Naghdi model on faceted triangulations.
By recovering directors numerically, we remove the requirement for an analytic or spline
parameterisation of the mid-surface, whereas partial selective reduced integration alleviates locking.
The residual and tangent of symbolic energies are generated in an automated fashion
building on Firedrake's and UFL's symbolic differentiation capabilities.
Additionally, UFL's and pyadjoint's algorithmic differentiation enables automated
first and second shape derivative evaluations that are verifiably correct via Taylor tests.
The forward model agrees with the Sze/Abaqus semi-cylinder reference.
A three-model mesh-refinement study indicates that geometry recovery and surface
faceting have small effects on the predicted deflection.
Numerical experiments show that 100 trust-region iterations suffice to reduce the squared-displacement
functional of a semi-cylinder by \semicylReductionPct\,\si{\percent} and the elastic strain energy
of a bracket by \bracketReductionPct\,\si{\percent}. The latter result is validated against
a COMSOL benchmark. Non-uniqueness of optimal shapes is demonstrated by optimising
initial designs subject to small
perturbations.
Comparisons with linearised designs distinguish optimisation performance from
the accuracy of the predicted nonlinear response.
\end{abstract}


\end{frontmatter}

\section*{Highlights}
\begin{itemize}\setlength{\itemsep}{0pt}
\item Nonlinear Naghdi shell on faceted meshes with a numerically recovered director
\item Geometry is the design variable: no analytic mid-surface chart required
\item Second-order-consistent shape derivatives from one recorded solve at the design load
\item Error decomposition: geometry recovery and faceting negligible against discretisation
\item Validated on the Sze/Abaqus forward and the COMSOL shell optimisation benchmark
\end{itemize}

\section{Introduction}
\label{sec:intro}

Thin-walled shells are ubiquitous in lightweight engineering because curved sheets of
material carry load predominantly through membrane action and so attain a high
stiffness-to-weight ratio. Because performance depends so strongly on curvature,
the \emph{shape} of a shell is a powerful design lever, and shape optimisation is
a natural route to lighter, stiffer components.

When the operating loads are large, the structural response of shells becomes
highly nonlinear because finite rotations, membrane--bending coupling
and load-dependent stiffening (or softening, up to snap-through) all matter.
In this regime, linearised models are not suitable since they may give an
inaccurate assessment of the shell response.
Therefore, shape optimisation loops must incorporate
geometrically nonlinear forward models, which significantly increases complexity.

Thanks to recent scientific software developments, tackling such complex optimisation problems
is now possible.
First, automated finite-element environments based on a symbolic weak-form language,
such as Firedrake \cite{rathgeber2016firedrake} and the Unified Form Language
\cite{alnaes2014unified}, generate the residual and the consistent tangent models
directly from energy functionals. Additionally, discrete adjoints can be derived
in an automated fashion through algorithmic differentiation \cite{mitusch2019pyadjoint}.
These libraries completely bypass the need for manual implementations, which
has historically made simulating nonlinear shell elements and their sensitivities laborious
and error-prone. A compelling demonstration of the extent of this simplification is demonstrated in
 the work on nonlinear Naghdi shells by Hale et~al.\ \cite{hale2018simple}.
 Second, shape-optimisation toolboxes such as
Fireshape \cite{paganini2021fireshape} provide a rigorous diffeomorphism-based
control space and Sobolev descent metrics that mitigate mesh-dependence and
promote mesh quality, so that the domain itself (instead of a low-dimensional parametrisation)
can be effectively treated as the optimisation variable.

Shape optimisation of shells has a long history \cite{ding1986shape,bletzinger1993form},
and geometrically nonlinear shells with adjoint sensitivities are part of it.
The contemporary literature is dominated by isogeometric analysis (IGA) approaches,
often based on the Kirchhoff--Love shell model
\cite{wall2008isogeometric,cho2009isogeometric,kiendl2014isogeometric,
hirschler2019embedded,zhao2024shape}. Applications range from
aerostructures to wind-turbine blades \cite{hirschler2019isogeometric}.
These approaches optimise shell shapes through their CAD parameters and
employ finite-dimensional optimisation algorithms. By definition, the number
of CAD parameters limits the dimension of the shape space available to the
optimiser. Enriching the shape space while retaining performance is not
straightforward due to a phenomenon known as \emph{mesh dependence}:
the number of optimisation steps to convergence may increase with the number
of dimensions considered \cite{SchHaFuPi17}. Possible approaches to mitigate
this undesirable effect are regularisation by multiresolution subdivision control meshes
\cite{bandara2018isogeometric} or
reducing the number of deformation variables, as in the IGA-based, open-source shape optimisation framework GOLDFISH \cite{benzaken2024goldfish}, which restricts control variables to free-form deformation parameters.
Antil et~al.\ \cite{antil2017shape} pursue a different approach.
Instead of using IGA, they manually reformulate the Naghdi shell model on a reference domain
through a midsurface chart. Shape optimization is then performed by directly optimizing this chart
(and also the shell thickness).
The shape derivatives of the resulting reference-domain operators are derived explicitly, and the chart
is approximated using piecewise-linear finite elements in the numerical implementation.

Instead, we are interested in an approach that does not require midsurface charts,
but operates directly on faceted triangulations, as these are commonly used in mechanical simulations. A notable previous
work in this direction is Firl and Bletzinger \cite{firl2012shape}, which considers a
continuum-based six-parameter Reissner--Mindlin shell element whose nodal director averages the
neighbouring element normals, so that no reference curvature is ever required.
Firl and Bletzinger apply optimisation directly to a node-based parametrisation of thin-walled structures,
assemble adjoint sensitivities from exact semi-analytical element derivatives,
employ convolution filters for regularisation,
and pursue a solution strategy that follows the load path of the initial design once and
then applies the design updates at fixed load.
Our approach is similar, but instead of relying on element-level semi-analytical formulas, we include
geometry operations within the differentiated program, thus completely bypassing
the need for manual derivation of derivative formulas. This change is significant
because it enables researchers to modify and compare models with ease,
more rapidly, and without concerns about the correctness of their implementation.
To demonstrate the validity of the approach, we consider the more general Naghdi shell
model, although the same approach can be applied to the Reissner--Mindlin model.
The contribution of
this paper is therefore neither a new shell theory nor a new optimiser, but
rather a shape optimisation framework built on a surface-tensor Naghdi model that is
chart-free and end-to-end differentiable.

This work is organised as follows. Sections~\ref{sec:model}--\ref{sec:fem}
describe the nonlinear Naghdi shell model formulation on \emph{discrete} (faceted) meshes
and its finite element implementation. Section~\ref{sec:shapeopt}
describes the formulation of a shape optimisation model for Naghdi shells
using the shape optimisation library Fireshape. Section~\ref{sec:results}
presents an extensive validation of the framework and includes:
validating the forward model against the Sze/Abaqus benchmark \cite{sze2004popular},
validating the shape optimisation method against the COMSOL ``Shape Optimization of a Shell'' benchmark
\cite{comsol_shell_shapeopt},
assessing the dependence of optimised designs on the initial geometry and mesh resolution,
and a comparison of linear versus nonlinear models.
The entire forward and adjoint model is expressed in fewer than a hundred lines
of Firedrake/UFL, which we reproduce in abridged form; the scripts that produce
every result of this paper are available as supplementary material.
\section{Nonlinear Naghdi shell model}
\label{sec:model}
We consider the five-parameter nonlinear Naghdi shell model. Both the Naghdi
and the Reissner--Mindlin models hinge on the same kinematic hypothesis of
first-order shear-deformation. Chapelle and Bathe \cite{chapelle2011finite}
show that the Naghdi model coincides with the shear-membrane-bending model,
and so both the Naghdi and the Reissner--Mindlin models encode the same physics
despite having different names.
What differs is their formulation. The continuum-based (degenerated-solid)
Reissner--Mindlin elements interpolate a three-dimensional body. They never
require the curvature of the mid-surface, but their sensitivities must be derived
element by element. Instead, the mid-surface formulation adopted here writes
the energy as a scalar functional of the displacement and of two director angles
through the surface tensors, namely the membrane strain, the change of curvature,
and the shear strain.
This form can be differentiated automatically using a symbolic finite-element language
like UFL
(also on standard Lagrange elements). This comes at the extra cost of
requiring the reference metric, director and curvature. In Section~\ref{sec:fem},
we describe how to recover this information from faceted meshes.
A third alternative is the Kirchhoff--Love model, which does not comprise director angles
but, in a standard conforming formulation, requires $C^1$-conforming discretisations
because its bending term depends on second derivatives of the displaced geometry.
Constructing finite elements that are $C^1$-conforming on triangulated meshes
is significantly more complex than the standard ($C^0$-conforming) Lagrangian counterparts,
and generally requires higher polynomial degrees. This explains why the Kirchhoff--Love model is often
discretised using an IGA method.
Alternatively, one can employ a $C^0$ interior penalty (C0IP)
formulation, as demonstrated for linear Koiter shells by Ferchichi~\cite{ferchichi2026interior}.
A detailed comparison of these three models within a shape optimisation framework is beyond the scope of this
work, and we simply consider the Naghdi shell model. However, we note that the shape-optimisation
methodology described in this work can be applied also to the Reissner--Mindlin and the Kirchhoff--Love models.

Next, we summarise the kinematics, the strain and stress measures,
and the resulting variational problem. For a detailed treatment of the Naghdi model we
refer the reader to \cite{naghdi1973theory,chapelle2011finite,hale2018simple}.
The reference configuration of the shell is the swept surface
\begin{equation}
\bm{p}_0(x_\alpha,x_3) = \bphi_0(x_\alpha) + x_3\,\bn_0(x_\alpha),
\qquad x_\alpha\in\omega\subset\Rbb^2,\quad x_3\in[-t/2,t/2],
\end{equation}
where $\bphi_0$ maps the parametric domain $\omega$ to the mid-surface and
$\bn_0$ is the unit normal. The deformed configuration is
\begin{equation}
\bm{p}(x_\alpha,x_3) = \bphi(x_\alpha) + x_3\,\bd(x_\alpha),
\qquad \bphi = \bphi_0 + \bu,
\end{equation}
where $\bu$ denotes the mid-surface displacement and $\bd$ the (unit, inextensible)
director. Parameterising $\bd$ by two angles
$\bbeta=(\beta_1,\beta_2)\in[-\tfrac{\pi}{2},\tfrac{\pi}{2}]\times(-\pi,\pi]$
as
\begin{equation}
\bd(\bbeta) = \bigl(\sin\beta_2\cos\beta_1,\,-\sin\beta_1,\,
\cos\beta_2\cos\beta_1\bigr)\,.
\label{eq:director}
\end{equation}
gives the five-parameter shell model. The kinematics of the model are illustrated in Figure~\ref{fig:kinematics}.

\begin{figure}[htb!]
\centering
\begin{tikzpicture}[>=Latex,scale=1.0]
  \draw[very thick,domain=-2:2,smooth,variable=\x]
        plot ({\x},{0.15*\x*\x-0.3});
  \node at (2.3,0.0) {$\bphi_0$};
  \coordinate (P0) at (-0.2,{0.15*0.04-0.3});
  \draw[->,shapeblue,very thick] (P0) -- ++(0.06,1.0) node[above left,yshift=-2pt]{$\bn_0$};
  \fill (P0) circle (1.2pt);
  \draw[very thick,red!70,domain=-2:2,smooth,variable=\x]
        plot ({\x+0.4},{0.45*\x+0.05*\x*\x+0.6});
  \node[red!70] at (3.35,1.45) {$\bphi=\bphi_0+\bu$};
  \coordinate (P) at ({-0.2+0.4},{0.45*(-0.2)+0.05*0.04+0.6});
  \draw[->,red,very thick] (P) -- ++(0.55,0.8) node[above]{$\bd(\bbeta)$};
  \fill (P) circle (1.2pt);
  \draw[->,gray,dashed,very thick] (P0) -- (P) node[midway,below right,yshift=3pt,xshift=4pt]{$\bu$};
\end{tikzpicture}
\caption{Kinematics of the nonlinear Naghdi shell: reference mid-surface
$\bphi_0$ with unit normal $\bn_0$, deformed mid-surface $\bphi$ with director
$\bd(\bbeta)$ from~\eqref{eq:director}, and the displacement $\bu$.}
\label{fig:kinematics}
\end{figure}

These kinematics fix the strain and stress measures. With
$\nabla\bphi=[\,\bphi_{,1}\;\bphi_{,2}\,]$ and reference metric
$\bm{a}_0=\nabla\bphi_0^{\!\top}\nabla\bphi_0$, the Lagrangian membrane, bending
and shear strains are
\begin{align}
\bm{e}(\bu)        &= \tfrac12\bigl(\nabla\bphi^{\!\top}\nabla\bphi-\bm{a}_0\bigr),\\
\bm{k}(\bu,\bbeta) &= -\tfrac12\bigl(\nabla\bphi^{\!\top}\nabla\bd
                       +\nabla\bd^{\!\top}\nabla\bphi\bigr)-\bm{b}_0,\\
\bm{\gamma}(\bu,\bbeta) &= \nabla\bphi^{\!\top}\bd - \nabla\bphi_0^{\!\top}\bd_0,
\end{align}
where $\bm{b}_0$ is the reference curvature. For a homogeneous
St.\,Venant--Kirchhoff material, the plane-stress elasticity tensor reads
\begin{equation}
A^{\alpha\beta\rho\sigma} =
\frac{2\lambda\mu}{\lambda+2\mu}\,a_0^{\alpha\beta}a_0^{\rho\sigma}
+\mu\bigl(a_0^{\alpha\rho}a_0^{\beta\sigma}+a_0^{\alpha\sigma}a_0^{\beta\rho}\bigr),
\end{equation}
giving the stress resultants
\begin{equation*}
\bm{N}=t\,\bA\!:\!\bm{e}\,,\quad
\bm{M}=\tfrac{t^3}{12}\bA\!:\!\bm{k}\,,
\quad\text{and}\quad
\bm{T}=t\,\mu\,\bm{a}_0^{-1}\bm{\gamma}\,,
\end{equation*}
and the energy densities
\begin{equation*}
\psi_m=\tfrac12\bm{N}\!:\!\bm{e}\,,\quad
\psi_b=\tfrac12\bm{M}\!:\!\bm{k}\,\quad
\text{and} \quad \psi_s=\tfrac12\bm{T}\!\cdot\!\bm{\gamma}\,.
\end{equation*}

Finally, equilibrium is derived variationally as the stationary configuration of the total potential
\begin{equation}\label{eqref:energy}
\Pi(\bu,\bbeta) = \int_\omega (\psi_m+\psi_b+\psi_s)\,\mathrm{d}x - W_{\mathrm{ext}},
\end{equation}
where $W_{\mathrm{ext}}$ is the external load work.
In Section~\ref{sec:results} we consider two load models.
For the semi-cylinder test case depicted in Figure~\ref{fig:semi_cylinder_setup},
the load is a  narrow Gaussian line load along that edge,
$W_{\mathrm{ext}}=\int_{\Gamma_{\mathrm{free}}} P\,g(x)\,u_z\,\mathrm{d}s$, where
$g$ is a unit-mass Gaussian of width $\sigma=\rho\pi/80$ in the circumferential
coordinate $x=\rho\sin\theta$, centred at the crown $\theta=0$.
This corresponds to a regularised central point load on the free edge
and avoids the otherwise resulting pointwise singularity (but converges to the point load as its width vanishes).
For the bracket test case depicted in Figure~\ref{fig:bracket-shape-sketch},
a total force $F$ is applied in the $y$-direction, distributed uniformly
over the loaded edge, that is,
$W_{\mathrm{ext}}=\int_{\Gamma_{\mathrm{load}}}(F/|\Gamma_{\mathrm{load}}|)\,
u_y\,\mathrm{d}s$. (Note that Figure~\ref{fig:semi_cylinder_setup} and Figure~\ref{fig:bracket-shape-sketch}
use different Cartesian axes.)

\section{Finite element implementation}
\label{sec:fem}

The distinguishing feature of the implementation used in this work is that it operates on
a \emph{discrete} (faceted) triangulation $\mathcal{T}_h$ of the mid-surface,
rather than on an exact analytic chart $\bphi_0$; every reference-geometry
quantity is therefore recovered from the mesh. On each flat triangle, the affine
geometric transformation that maps the finite element reference triangle
onto the physical triangle has a constant Jacobian matrix $J_0$. Specifically,
$J_0$ is a $3\times2$ matrix whose two columns are the in-plane edge tangents.
In this case, the reference metric (i.e., the first
fundamental form) is the constant $\bm{a}_0=J_0^{\!\top}J_0$. The director $\bd(\bbeta)$ is
built from the surface normal. Each triangle carries a single constant unit
normal, its \emph{cell normal} $\bm{\nu}$ (UFL's \texttt{CellNormal}),
but these per-facet normals are discontinuous across shared edges. To recover
a continuous normal field, we compute the $L^2$-projection of
$\bm{\nu}$ onto the continuous piecewise-linear Lagrangian space $[\mathrm{CG}_1]^3$,
which averages the adjacent facet normals at each vertex, and then renormalise
at the mesh nodes. Formally, we compute
\begin{equation}
\bn_0 = \frac{\mathcal{P}_h\,\bm{\nu}}{\lvert\mathcal{P}_h\,\bm{\nu}\rvert}
\quad\text{and}\quad
\bbeta_0 = \bigl(\,\mathrm{atan2}(-n_{0,2},\sqrt{n_{0,1}^2+n_{0,3}^2})\,,
\mathrm{atan2}(n_{0,1},n_{0,3})\,\bigr),
\end{equation}
where $\mathcal{P}_h$ denotes the $L^2$-projection onto $[\mathrm{CG}_1]^3$ and
$\bbeta_0=(\beta_{0,1},\beta_{0,2})$ are the two director angles read off the
recovered normal $\bn_0$. The reference director $\bd_0=\bd(\bbeta_0)$ can be
computed from these angles using formula~\eqref{eq:director}, and the reference curvature
(i.e., the second fundamental form) is recovered consistently from its facet gradient,
\begin{equation}
\bm{b}_0 = -\tfrac12\bigl(J_0^{\!\top}\nabla\bd_0
           + \nabla\bd_0^{\!\top}J_0\bigr)\,.
\label{eq:b0-recovered}
\end{equation}
Following this procedure, the metric $\bm{a}_0$, the curvature $\bm{b}_0$, and the director $\bd_0$
are all obtained from the mesh alone.
This recovered-geometry route is what makes the model compatible with shape
optimisation: although every design update produces a new mesh with no closed-form
parameterisation, these reference quantities remain well defined. To guarantee
a coherent director sign across the whole triangulation, we simply fix (once) a
consistent global orientation of the cell normal field $\bm{\nu}$.
Although the recovery procedure we just described is not exact,
the resulting geometric error is consistent and vanishes under mesh refinement
\cite{hale2018simple}, as evidenced numerically in Section~\ref{sec:results-validation}.

On this recovered geometry, the displacements are discretised with bubble-enriched
quadratic Lagrangian finite elements (that is, $\bu_h\in[\mathrm{CG}_2\oplus B_3]^3$)
and rotations with quadratic Lagrangian finite elements (that is, $\bbeta_h\in[\mathrm{CG}_2]^2$).
Membrane and shear locking are alleviated with a
partial selective reduced integration (PSRI) scheme: the membrane and shear
energies are split into full- and reduced-quadrature parts weighted by a cell
factor $\alpha$,
\begin{equation}
\Pi_h = \Pi_b
 + \alpha(\Pi_m+\Pi_s)\big|_{\text{full}}
 + (1-\alpha)(\Pi_m+\Pi_s)\big|_{\text{red}},
\qquad \alpha = t^2/h_{\mathrm{eff}}^2,
\label{eq:psri}
\end{equation}
with full quadrature of degree~4 and reduced quadrature of degree~2.
Specifically, we use
$h_{\mathrm{eff},K}^2=c_K|\bm j_1\times\bm j_2|$, where
$c_K=h_K(0)^2/|\bm j_1(0)\times\bm j_2(0)|$ is fitted on the initial mesh and
held fixed. Thus $h_{\mathrm{eff}}$ equals the cell diameter initially and
thereafter varies differentiably with the facet area; it is not the current
cell diameter under arbitrary distortion. The scaling $\alpha\simeq t^2/h^2$
follows \cite{hale2018simple,lovadina2005energy} and balances membrane/shear
and bending terms.

Finally, to compute the equilibrium point of \eqref{eqref:energy}, we use load continuation:
we consider a sequence of incremental loadings $P\in\{0,\dots,P_{\max}\}$ and sequentially
apply Newton's method (PETSc SNES, \texttt{newtonls} with an $\ell^2$ line search and a MUMPS LU factorisation)
using each converged state as the initial guess for the next load.
The whole discretisation is expressed compactly in Firedrake.
Listing~\ref{lst:forward} shows the core of the forward model. Notably,
the tangent and its Jacobian are obtained from the energy
by automatic differentiation
(command \texttt{derivative}); no hand-coded formula is required.

\begin{lstlisting}[style=python,caption={Core of the nonlinear Naghdi forward
model in Firedrake (abridged).},label={lst:forward},
moredelim={[is][\color{red}]{@red@}{@endred@}}]
# Five-parameter director from two angles
def director(beta):
    return as_vector([sin(beta[1])*cos(beta[0]),
                      -sin(beta[0]),
                      cos(beta[1])*cos(beta[0])])

# Reference geometry recovered from the faceted mesh
J0 = Jacobian(mesh); a0 = J0.T*J0; a0_contra = inv(a0)
n0 = project(CellNormal(mesh), V_normal)                       # recovered normal
n0 = Function(V_normal).interpolate(n0 / sqrt(inner(n0, n0)))
beta0 = interpolate(as_vector([atan2(-n0[1], sqrt(n0[0]**2 + n0[2]**2)),
                               atan2(n0[0], n0[2])]), V_beta)  # initial director angles
d0 = director(beta0)

# Strain measures (membrane e, bending k, shear gamma)
F = dot(grad(u_), J0) + J0
d = director(beta_ + beta0)
e = 0.5*(F.T*F - a0)
k = -0.5*(F.T*dot(grad(d), J0) + dot(grad(d), J0).T*F) - b0
gamma = F.T*d - J0.T*d0

# PSRI: split membrane/shear into full (deg 4) and reduced (deg 2) parts
area_scale = sqrt(dot(cross(J0[:, 0], J0[:, 1]),
                      cross(J0[:, 0], J0[:, 1])))
# c_K = initial_cell_diameter**2 / initial_area_scale, held fixed
alpha = t**2 / (c_K * area_scale)
Pi = (psi_b*dx_full
      + alpha*(psi_m + psi_s)*dx_full
      + (1.0 - alpha)*(psi_m + psi_s)*dx_red)

R = derivative(Pi, q_, q_t) + W_ext  # residual (auto-differentiated)
J = derivative(R, q_, q)             # consistent tangent
solve(R == 0, q_, bcs=bcs, J=J, solver_parameters={...})  # solver parameters abridged
\end{lstlisting}

\section{Shape optimisation implementation}
\label{sec:shapeopt}

Optimising the shell design to increase performance is a PDE-constrained shape-optimisation problem:
the optimisation variable is the domain itself, and every objective evaluation requires a full nonlinear shell
solve. Following \cite{paganini2021fireshape,delfour2011shapes}, we describe domains as
images of bi-Lipschitz transformations of the form
$\bT:\Rbb^d\to\Rbb^d$ applied to an initial domain $\Omega$,
so that the set of admissible domains is
\begin{equation}\label{eq:admissible_domains}
\mathcal{Q} = \{\,\bT(\Omega) : \bT \text{ bi-Lipschitz}\,\}.
\end{equation}

To model the mechanical performance of a shell design we consider different objectives
for different test cases to better reflect the tailored measures of stiffness.
For the semi-cylinder test case, we minimise the \emph{squared-displacement
functional} $J_{\mathrm{sc}}\coloneqq\int_\omega \lvert\bu\rvert^2\,\mathrm{d}x$,
which is a direct measure of the kinematic response under load.
For the bracket test case, we minimise the total elastic strain energy (i.e., the compliance)
$J_{\mathrm{br}}\coloneqq\int_\omega(\psi_m+\psi_b+\psi_s)\,\mathrm{d}x$, which was the
objective used in the COMSOL benchmark. Note that
the same adjoint pipeline can handle either objective.

To discourage pointwise shape displacements beyond
a prescribed budget $\delta_{\max}$,
we introduce the penalty term
\begin{equation}
P_\delta(\bT) \coloneqq \tfrac12\int_\omega
\Bigl[\bigl(\lvert\bT-\bI\rvert-\delta_{\max}\bigr)_{+}\Bigr]^{2}
\,\mathrm{d}x,
\qquad (s)_{+}\coloneqq\max(s,0)\,,
\label{eq:shape-penalty}
\end{equation}
where $\bI$ denotes the identity operator $\bI(x)=x$.
In the implementation,
the norm and the positive part $(\,)_+$ are smoothed with a small parameter $\varepsilon$
(e.g.\ $\lvert\cdot\rvert\to\sqrt{\lvert\cdot\rvert^{2}+\varepsilon^{2}}$) and
we replace $\delta_{\max}$ with $\delta_{\max}+\varepsilon$, so that $P$ is everywhere
weakly differentiable and
$P_\delta(\bI)=0$. A displacement budget bounds where the surface may go
but not how much it may fold within that bound. Therefore, we add an additional
soft budget on the relative change of the area $A(\bT)\coloneqq\int_{\omega(\bT)}\mathrm{d}x$
of the design surface. Specifically, we consider the penalty term
\begin{equation}
P_{A}(\bT) \coloneqq \tfrac12\Bigl[\bigl(A(\bT)/A_0 - 1 - a_{\max}\bigr)_{+}\Bigr]^{2},
\label{eq:area-penalty}
\end{equation}
where $A_0$ is the initial area and $a_{\max}$ is the admissible relative increase.
In the implementation, Formula~\eqref{eq:area-penalty} is smoothed in a similar
fashion to Formula~\eqref{eq:shape-penalty}.
Its shape derivative follows from applying the chain rule on
the first variation of the area, expressed as a UFL form on the moved mesh.
The Hessian
action introduces only a rank-one term and therefore requires just two scalar assemblies per evaluation,
which is negligible.
The final objective is a weighted sum $J$ of the original objective $J_{\mathrm{sc}}$
or $J_{\mathrm{br}}$ and (a smoothed version of) the penalisation terms $P_\delta$ and $P_A$.
In the semi-cylinder test case, we also consider an additional penalty on the tangential
control slope near its mirror plane,
as described in Section~\ref{sec:results-validation}. Note that the penalty terms are
so-called ``soft penalties'':  although they heavily penalise constraint violation,
they do not guarantee that nominal budgets are satisfied exactly.

The optimiser updates are driven by the shape derivative $dJ$ of the
objective $J$, computed by algorithmic differentiation with UFL and
pyadjoint \cite{mitusch2019pyadjoint,ham2019automated}. The objective depends on
the shape only through the equilibrium state at the full load, so it is that
state, and not the path that leads to it, which has to be differentiated. The
forward program is therefore recorded once, at the reference shape, as a single
Newton solve at the full load preceded by the geometry recovery of
Section~\ref{sec:fem} and followed by the assembly of the objective. For every
trial shape the recorded program is replayed with the new coordinates: the
geometry is recovered, the Newton solve is started from the equilibrium state
on the previous geometry and, after convergence, the objective is evaluated directly
while its derivative is obtained from one shell-state and geometry-recovery adjoint solve.
The load continuation procedure described in Section~\ref{sec:fem} is used only in two instances:
to generate the equilibrium state at the initial shape, and to generate a better initial guess
for the Newton solve if it fails to converge starting from the equilibrium computed on
the previous shape iterate.

To guarantee that only feasible designs are considered, we monitor the
surface's element-wise area, which changes sign if an element is inverted.
Since the shell is a two-dimensional manifold
embedded in $\Rbb^3$, this feasibility measure is computed using the formula
\begin{equation*}
\sqrt{\det(J_0^{\!\top}J_0)}=\lvert\bm{j}_{0,1}\times\bm{j}_{0,2}\rvert\,,
\end{equation*}
where $\bm{j}_{0,1}$ and $\bm{j}_{0,2}$ are the columns of the Jacobian $J_0$
of the moved element map; cf. Section~\ref{sec:fem}.
A trial shape is immediately rejected if it fails this feasibility test.
The mesh quality of optimised designs is verified in Section~\ref{sec:results-verification}.
Note that this local feasibility test cannot detect nonlocal self-contact nor impose
triangle-quality bounds.

To discretise the set of admissible domains $\mathcal{Q}$ defined in \eqref{eq:admissible_domains},
we consider transformations of
the form $\bm{T}=\bI+\bq$, where $\bq$ is a boundary-preserving
vector field that displaces the mesh interior. This choice of transformations
enables ``discretising'' $Q$ by choosing to which finite element space $\bq$ should belong.
Because we are primarily interested in optimal shell designs that exhibit
large-scale shape changes instead of small-amplitude
high-frequency minima, we can select a relatively lower-dimensional
shell parametrisation, which is refined for the accurate computation
of the mechanical shell response. This approach is often called a \emph{multilevel}
(or geometric multigrid) discretisation \cite{paganini2021fireshape,schulz2016efficient}.
We can summarise its implementation as follows:
first, we construct a mesh hierarchy via uniform refinement of a coarse triangulation
$\mathcal{T}_0$ of the mid-surface, $\mathcal{T}_0 \to \cdots \to \mathcal{T}_L$;
then, we construct nested vector coordinate spaces $V_0\subset\cdots\subset V_L$.
The \emph{control} variable $\bq$ lives on the coarsest level, that is, $\bq\in V_0$,
whereas the shell state is computed using the
\emph{finest} level $\mathcal{T}_L$. To transfer control updates to the finest mesh $\mathcal{T}_L$,
we employ the  \emph{prolongation} (the composition of inter-grid transfers) operator
$E:V_0\to V_L$ and deform the finest mesh using the geometric transformation
\begin{equation}
\bm{T} = \bI + E(\bm{q})\,.
\end{equation}
The shape derivative of ${J}$ is first evaluated along mesh deformations of
the fine mesh using the adjoint method. Then, applying the transpose of the
prolongation $E$ (the standard multigrid restriction), we obtain its
\emph{restriction} onto the coarse control space.

Finally, we should endow $Q$ with an \emph{inner product}
$a(\cdot\,,\cdot)$ to compute the shape gradient $\bm{\nabla J}\in V_0$, which
informs how the optimiser updates $\bq$ and is the solution to
\begin{equation}
a(\bm{\nabla{J}}, \bw) = d{J}(\bw)\quad\text{for all } \bw\in V_0\,.
\end{equation}
Since we are interested in smooth, large-scale shape changes, we consider the
\emph{smoothing} (Helmholtz) metric,
\begin{equation}
a(\bg,\bw)=\int_{\Omega}\big(\bg\cdot \bw + \ell^{2}\,\nabla \bg:\nabla \bw\big)\,\mathrm{d}x ,
\label{eq:helmholtz-metric}
\end{equation}
which acts as a Helmholtz filter of feature size
$\ell$, closely related to the Helmholtz/PDE filters and vertex-morphing
regularisation used in node-based structural shape optimisation
\cite{lazarov2011helmholtz,hojjat2014vertex}: a larger length scale $\ell$
penalises short-wavelength control variation more strongly and so sets the
wavelength of the folds. Because $\dim V_0\ll\dim V_L$,
the design space automatically excludes high-frequency designs, while the fine
state mesh enables accurate computation of the energy and its gradient using the
faceted forward model described in Section~\ref{sec:fem}.
In Fireshape this approach is implemented via the \texttt{FeMultiGridControlSpace} and a
Helmholtz \texttt{UflInnerProduct} (Listing~\ref{lst:fireshape}).
\begin{lstlisting}[style=python,caption={Coupling the shell to Fireshape and ROL
with a multilevel control and a Helmholtz smoothing metric (abridged).},
label={lst:fireshape}]
class ShellShapeObjective(SingleTapeShellObjective):
    def objective_value(self):
        self.update_geometry_coefficients()    # recompute n0, beta0, alpha
        self.P_load.assign(self.max_load)
        self.solver.solve()                   # one full-load solve on tape
        return assemble(self.psi_total*dx) / self.reference_objective
# Untaped continuation initialises the state. If replay fails, continuation
# supplies a new initial guess and the same full-load tape is replayed.

class HelmholtzInnerProduct(UflInnerProduct):  # smoothing Riesz map (M + l^2 K)
    def get_weak_form(self, V):
        u, v = TrialFunction(V), TestFunction(V)
        return (self.l**2*inner(grad(u), grad(v)) + inner(u, v))*dx

# Multilevel control on a mesh hierarchy: coarse control, fine state
mh = MeshHierarchy(base_mesh, refine)
Q = FeMultiGridControlSpace(mh, coarse_control=True)
inner = HelmholtzInnerProduct(Q, length_scale=..., fixed_bids=FIXED_TAGS,
                              extra_bcs=symmetry_bcs, direct_solve=True)
q = ControlVector(Q, inner)
J = (ShellShapeObjective(Q, ...)
     + ShapeDisplacementPenalty(Q, limit=..., smoothing=..., scale=...)
     + SurfaceAreaPenalty(Q, max_area_change=..., smoothing=..., scale=...))
# The semi-cylinder also adds SymmetryDerivativePenalty.

# Trust-region BFGS via ROL
params["General"]["Secant"]["Use as Hessian"] = True
params["General"]["Secant"]["Maximum Storage"] = 10
solver = ROL.OptimizationSolver(ROL.OptimizationProblem(J, q), params)
solver.solve()
\end{lstlisting}

The resulting optimisation problem is solved with the optimisation library ROL
\cite{rol2014}, which already includes a large variety of optimisation algorithms.
Specifically, we employ ROL's trust-region method with a Steihaug--Toint truncated
conjugate-gradient (CG) solver for the trust-region subproblem. We explicitly set
\texttt{Use as Hessian} to \texttt{True} and retain ten L-BFGS secant pairs to construct
the quadratic model. This choice is motivated by our previous experience with PDE-constrained
shape optimisation and provides a practical balance between robustness and computational cost,
rather than being
an intrinsic requirement of the proposed framework. In particular, the trust-region strategy controls
the magnitude of shape updates according to the agreement between predicted and actual objective
reductions. This is useful for geometrically nonlinear shell problems, where large shape updates may
cause mesh deterioration or failure of the nonlinear state solve.
Additionally, the L-BFGS approximation avoids the additional PDE solves required to evaluate exact second-order
derivatives, while its storage and application costs scale linearly with the number of shape-design variables
for a fixed number of secant pairs. The truncated-CG solver further allows the trust-region subproblem to
be solved without assembling a dense Hessian matrix.
Figure~\ref{fig:optloop} shows the overall optimisation loop.

\begin{figure}[htb!]
\centering
\begin{tikzpicture}[
  >=Latex, node distance=7mm and 13mm,
  box/.style={draw,rounded corners,align=center,minimum height=11mm,
              minimum width=27mm,font=\footnotesize},
  every edge/.style={draw,->}]
  \node[box] (rol) {ROL trust region\\(BFGS, TCG)};
  \node[box,right=of rol] (ctrl) {Control update\\$\bm{q}\to\bT(\Omega)$};
  \node[box,right=of ctrl] (mesh) {Move mesh};
  \node[box,below=of mesh] (fwd) {Forward solve};
  \node[box,left=of fwd] (obj) {Objective\\${J}(\bm{q})$};
  \node[box,left=of obj] (adj) {Adjoint\\$\bm{\nabla{J}}$};
  \draw (rol) edge (ctrl);
  \draw (ctrl) edge (mesh);
  \draw (mesh) edge (fwd);
  \draw (fwd) edge (obj);
  \draw (obj) edge (adj);
  \draw (adj) edge (rol);
\end{tikzpicture}
\caption{PDE-constrained shape-optimisation loop coupling ROL, Fireshape
(control space and mesh motion) and the Firedrake nonlinear Naghdi forward and
adjoint solves.}
\label{fig:optloop}
\end{figure}

Alternative optimisation algorithms, including line-search L-BFGS, nonlinear conjugate-gradient
and Newton-type methods, can also be employed within the ROL framework. SQP methods may be
appropriate when treating explicit design constraints, but here the only constraints beyond the
state equation are the soft penalties~\eqref{eq:shape-penalty} and~\eqref{eq:area-penalty}.
Importantly, ROL employs a
function-space-oriented optimisation based on the control-space metric, duality pairings and derivative
actions, which are appropriately set by Fireshape. Consequently, the choice of optimisation algorithm
is not dictated solely by the number of discrete shape-design variables, although discretisation, conditioning
and PDE-solution costs continue to influence computational performance.
A systematic comparison of optimisation algorithms could be of interest, but is beyond the scope of this work.
Nevertheless, in Section~\ref{sec:results-verification}, we compare the adopted L-BFGS approximation with the exact Hessian of the reduced objective, keeping the trust-region algorithm unchanged.

\section{Numerical experiments}
\label{sec:results}

We consider two test cases: the clamped semi-cylindrical shell used in the Sze
benchmark \cite{sze2004popular}, which serves both to validate the forward model
against the Abaqus S4R reference and as a shape-optimisation showcase, and
the sheet-metal bracket used in the COMSOL ``Shape Optimization of a Shell''
benchmark \cite{comsol_shell_shapeopt}, which validates the shape optimisation methodology.
For each test case, we first describe the geometry, boundary conditions, loading,
design constraints, meshes and optimiser settings before discussing the numerical results.
All computations use the discretisation described in Section~\ref{sec:fem}
(with relative and absolute residual tolerances set to $10^{-8}$ and $10^{-9}$
as stopping criteria for the Newton solver), and the
optimiser of Section~\ref{sec:shapeopt} with initial trust-region radius
$0.01$, radius growing rate $2$, shrinking rates $0.25$ and $0.5$, step-acceptance
ratio $10^{-3}$, gradient tolerance $10^{-3}$, step tolerance $10^{-7}$ and an
iteration limit of $100$. The soft shape-displacement penalty~\eqref{eq:shape-penalty}
and the area penalty~\eqref{eq:area-penalty} are
included in the objective $J$ using the multiplicative constant factors $10^{5}$
and $10^4$, respectively, and a relative smoothing of $10^{-3}$. The design reported for each
example is the final recorded iterate. Mesh quality is measured using the area penalty
(Section~\ref{sec:results-verification}).
The baseline runs were performed with eight MPI
ranks on one node with a 64-core AMD Ryzen Threadripper PRO 9985WX processor
and $124$\,GiB of memory.

\subsection{Clamped semi-cylindrical shell}
\label{sec:results-validation}

\paragraph{Problem definition}
We consider a semi-cylindrical (half-cylinder) shell whose mid-surface
$$\{(\rho\sin\theta,\,y,\,\rho\cos\theta):\theta\in[-\tfrac{\pi}{2},\tfrac{\pi}{2}],
\,y\in[0,L]\}$$ has radius $\rho$, axial length $L$ and uniform thickness $t$.
Geometry and material values are collected in Table~\ref{tab:semicyl-params}.
The setup is displayed in Figure~\ref{fig:semi_cylinder_setup}(a).
The shell is fully clamped ($\bu=\bm{0}$, $\bbeta=\bm{0}$) at the curved edge
located at $y=0$. The opposite (curved) edge, located at $y=L$, is left free.
The two straight lateral edges
$\theta=\pm\tfrac{\pi}{2}$ carry the mirror-symmetry conditions $u_z=0$ and
$\beta_2=0$ that let this half model represent a full cylinder. A point load
$P$, directed downward and raised by load continuation to a maximum of
\SI{2}{\kilo\newton}, is applied at the crown $\theta=0$ of the free edge.
Since this model is symmetric with respect to the plane $x=0$, for shape
optimisation we consider the half model $\theta\in[0,\pi/2]$
(complemented with the mirror-plane conditions $u_x(x=0)=0$ and $\beta_2(x=0)=0$
on the tangential shape displacement) and reflect results for visualisation.
We also include a soft penalty on the in-plane slope of the
control at the mirror plane to encourage a smooth reflected design across
the crown.
For shape optimisation, shape modifications are applied only to the interior of the mid-surface,
and not on its boundary. Shape modifications are also limited to
$\delta_{\max}=\SI{\dispBudgetSemicyl}{\meter}$, and the area increase is limited to
$a_{\max}=+\areaBudgetSemicyl\,\si{\percent}$. The control mesh is the $1024$-triangle half mesh,
and it is uniformly refined
once to the $4096$-triangle state mesh to compute shell displacements (which
corresponds to $54{,}213$ state degrees of
freedom). The objective to be minimised is the squared-displacement functional
$J_{\mathrm{sc}}$ of Section~\ref{sec:shapeopt} scaled by its value on the
reference cylinder. The Helmholtz length scale is $\ell=\SI{1.0}{\meter}$. The
load continuation that initialises the tape and recovers failed replays uses
$40$ increments.

\begin{figure}[htb!]
\centering
\input{semicyl_setup_sketch.tex}
\caption{Reference geometry, boundary conditions, loads, and fixed optimisation boundaries (in blue)
used in the semi-cylinder test case.}
\label{fig:semi_cylinder_setup}
\end{figure}

\begin{table}[t]
\centering\small
\caption{Clamped semi-cylindrical shell benchmark: geometry, material and
conditions.}
\label{tab:semicyl-params}

\begin{tabular}{@{}lll@{}}
\toprule
\multicolumn{3}{@{}l}{\textbf{Clamped semi-cylindrical shell}}\\
\midrule

\textit{Geometry}
&
\textit{Material}
&
\textit{Conditions}
\\[2pt]

\begin{tabular}[t]{@{}ll@{}}
Radius $\rho$
  & \SI{1.016}{\meter} \\
Axial length $L$
  & \SI{3.048}{\meter} \\
Thickness $t$
  & \SI{30}{\milli\meter}
\end{tabular}
&
\begin{tabular}[t]{@{}ll@{}}
Young's modulus $E$
  & \SI{20.685}{\mega\pascal} \\
Poisson ratio $\nu$
  & \num{0.3}
\end{tabular}
&
\begin{tabular}[t]{@{}ll@{}}
Applied load $P$
  & \SI{2}{\kilo\newton} \\
Shape budget $\delta_{\max}$
  & \SI{\dispBudgetSemicyl}{\meter} \\
Area budget $a_{\max}$
  & $+\areaBudgetSemicyl\,\si{\percent}$ \\
Helmholtz scale $\ell$
  & \SI{1.0}{\meter} \\
Control/state mesh
  & $1{,}024/4{,}096$ triangles
\end{tabular}
\\

\bottomrule
\end{tabular}
\end{table}

\begin{figure}[htb!]
\centering

\begin{subfigure}[c]{0.50\linewidth}
  \centering
  \includegraphics[width=\linewidth]{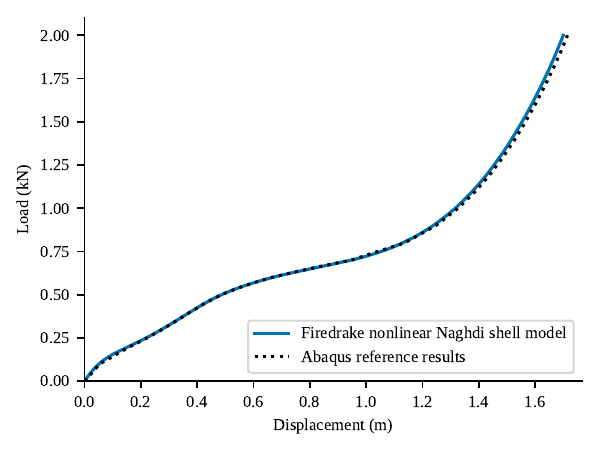}
  \caption{}
\end{subfigure}
\hspace{0.2cm}
%
\begin{minipage}[c]{0.36\linewidth}
  \centering

  \begin{subfigure}[c]{0.68\linewidth}
    \centering
    \includegraphics[width=\linewidth]{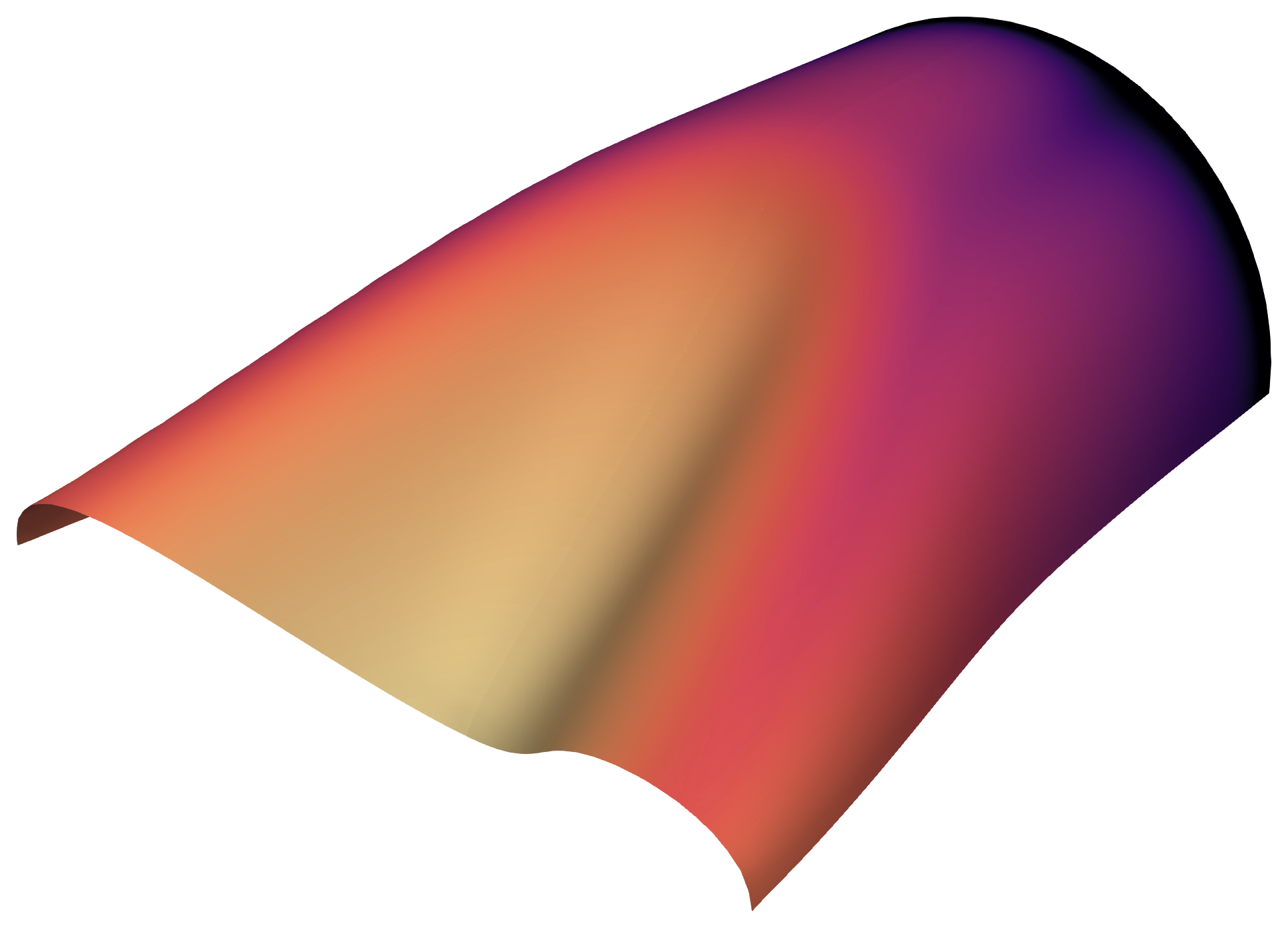}
    \caption{}
  \end{subfigure}
  \hfill
  \begin{minipage}[c]{0.18\linewidth}
  \end{minipage}

  \vspace{0.5em}


  \begin{subfigure}[c]{0.68\linewidth}
    \centering
    \includegraphics[width=\linewidth]{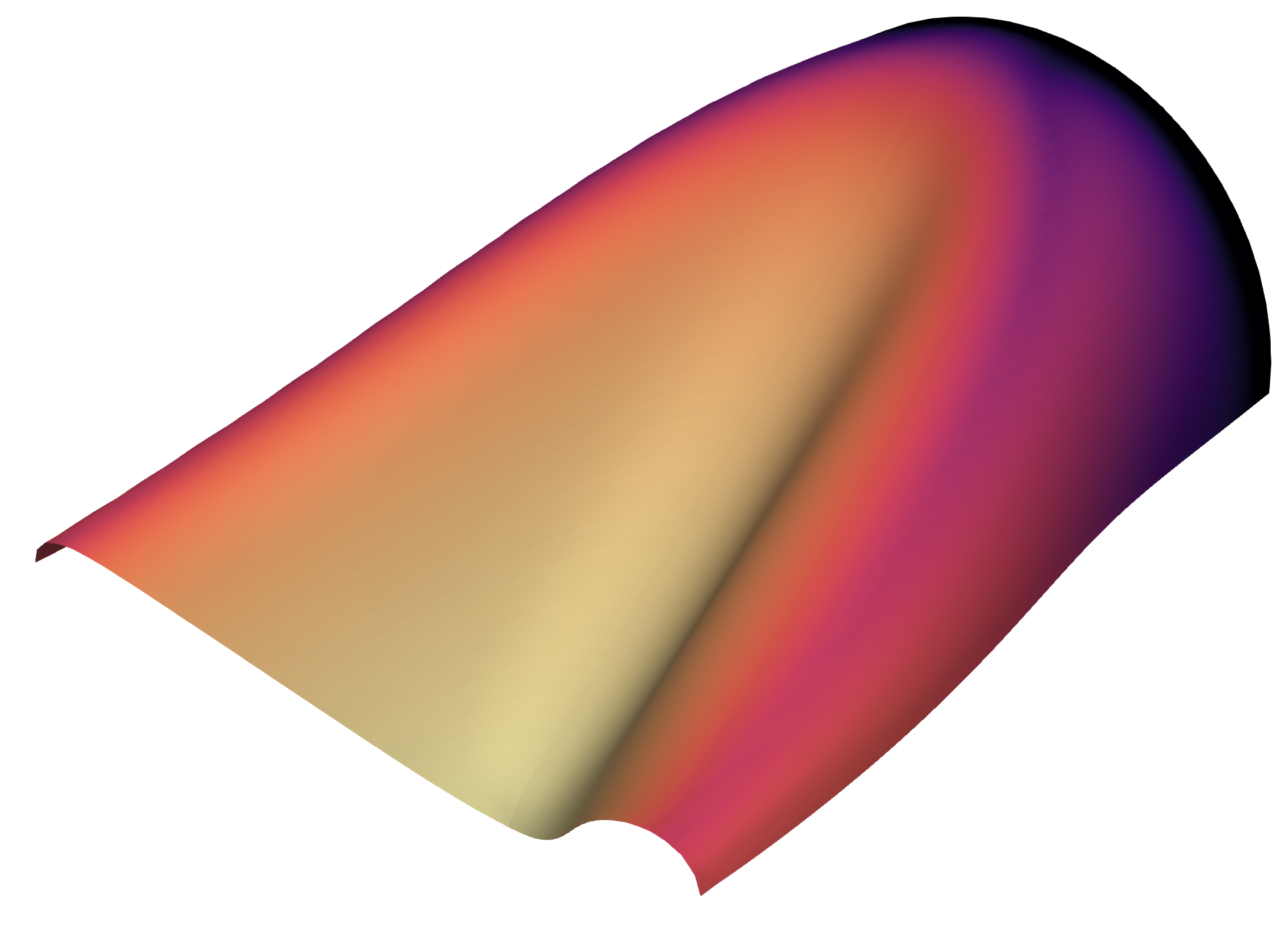}
    \caption{}
  \end{subfigure}
  \hfill
  \begin{minipage}[c]{0.18\linewidth}
    \centering
    \begin{tikzpicture}
      \begin{axis}[
          hide axis,
          scale only axis,
          height=3.1cm,
          width=1pt,
          colormap name=magma,
          colorbar,
          point meta min=-3,
          point meta max=0.30103,
          colorbar style={
            width=0.26cm,
            height=3.1cm,
            ytick={-2,-1,0},
            yticklabels={$0.01$,$0.1$,$1$},
            yticklabel pos=right,
            ytick pos=right,
            ylabel={Displacement (m)},
            ylabel near ticks,
            tick label style={font=\footnotesize},
            label style={font=\footnotesize},
            axis line style={line width=0.5pt, black},
            tick align=outside,
          }]
        \addplot[draw=none] coordinates {(0,0) (1,1)};
      \end{axis}
    \end{tikzpicture}
  \end{minipage}
  \hfill

\end{minipage}

\caption{Clamped semi-cylindrical shell under point load: (a) load--deflection of
the downward displacement at the load point, Firedrake nonlinear Naghdi model
vs.\ the Sze/Abaqus S4R reference \cite{sze2004popular},
(b-c) illustrative forward responses at \SI{1}{\kilo\newton} and \SI{2}{\kilo\newton}, respectively,
with the mesh warped by displacement and coloured by its magnitude.}
\label{fig:validation-deformed}
\end{figure}

\paragraph{Forward validation}
We compute the load--deflection curve of the shell and compare it to the
data from the Sze/Abaqus S4R benchmark \cite{sze2004popular}.
The results shown in Figure~\ref{fig:validation-deformed}(a) highlight an excellent match between
the two lines across the full load range, including the stiffening beyond
$\SI{750}{\newton}$. The agreement is also quantitative: the root-mean-square
error against the reference points is \fwdRMSE\,\si{\meter}
(\fwdRelErr\,\si{\percent} of the peak deflection), with a maximum error of
\fwdMaxErr\,\si{\meter}. Figure~\ref{fig:validation-deformed}(b) and (c) show
the deformed shell configurations at \SI{1}{} and \SI{2}{\kilo\newton}, respectively.

\begin{figure}[htb!]
\centering
\includegraphics[width=0.5\linewidth]{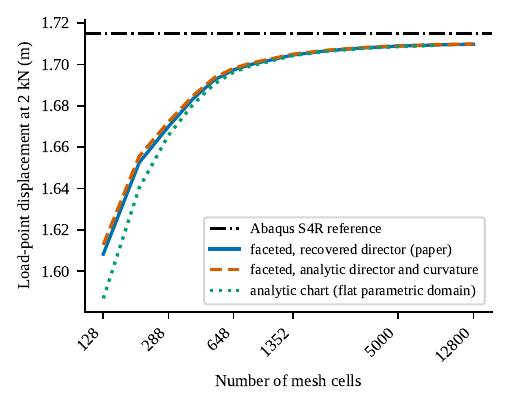}
\caption{Load-point deflection of the semi-cylinder at \SI{2}{\kilo\newton}
versus the number of mesh cells for the faceted model with recovered director
(Section~\ref{sec:fem}), the same faceted meshes with the analytic director and
curvature, and the chart-based model on the flat parametric domain; the dotted
line is the Sze/Abaqus S4R reference.}
\label{fig:forward-geometry}
\end{figure}

\paragraph{Discretisation, geometry and model errors}
To assess (and distinguish between) the discretisation error, faceted geometry
effects, and competing shell models, we compute the load-point displacement (m)
at \SI{2}{\kilo\newton} on
a sequence of meshes (ranging from $128$ to $12{,}800$ triangles)
using three forward models: the model of Section~\ref{sec:fem}, in which the director and curvature
are recovered from the triangulation; the same faceted mesh with the analytic
cylinder director and curvature evaluated at the quadrature points; and the
chart-based model of Hale et~al.\ \cite{hale2018simple}, in which the flat
parametric domain is mapped to the surface by the analytic function.
The results are portrayed in Figure~\ref{fig:forward-geometry}.
We observe that the recovery error is negligible: at the
\SI{800}{}-triangle validation mesh the recovered and analytic geometries give
load-point deflections at \SI{2}{\kilo\newton} of \SI{1.6999}{} and
\SI{1.7006}{\meter} (\SI{0.04}{\percent}), and at $12{,}800$ cells
\SI{1.70977}{\meter} and \SI{1.70985}{\meter} (\SI{0.005}{\percent}). The faceting of
the metric is not visible either: the chart-based model converges to
\SI{1.7092}{\meter} at $8192$ triangles against \SI{1.7095}{\meter} for the faceted
models. All three models approach the displacement \SI{1.710}{\meter}, which
is \SI{0.3}{\percent} below
the S4R value of \SI{1.715}{\meter},
and the root-mean-square errors (compared to the Abaqus' load-displacement
curve over the load path) are \SI{0.31}{\percent} (faceted) and \SI{0.28}{\percent} (chart) of
the maximal deflection. The validation mesh is within
\SI{0.5}{\percent} of the converged value. At coarse resolution the faceted
models are in fact closer to the converged value than the chart-based one
(\SI{1.608}{} and \SI{1.613}{\meter} against \SI{1.587}{\meter} at $128$ cells),
because their vertices lie on the cylinder while the chart is interpolated
quadratically. The point-load regularisation is likewise benign: halving the
Gaussian half-width (from $\rho\pi/80$ to $\rho\pi/160$) changes the
\SI{2}{\kilo\newton} deflection by \SI{0.21}{\percent}, and a fourfold range of
widths spans under \SI{0.4}{\percent}.

\begin{figure}[htb!!]
\centering


\begin{subfigure}[c]{0.31\linewidth}
  \includegraphics[width=\linewidth]{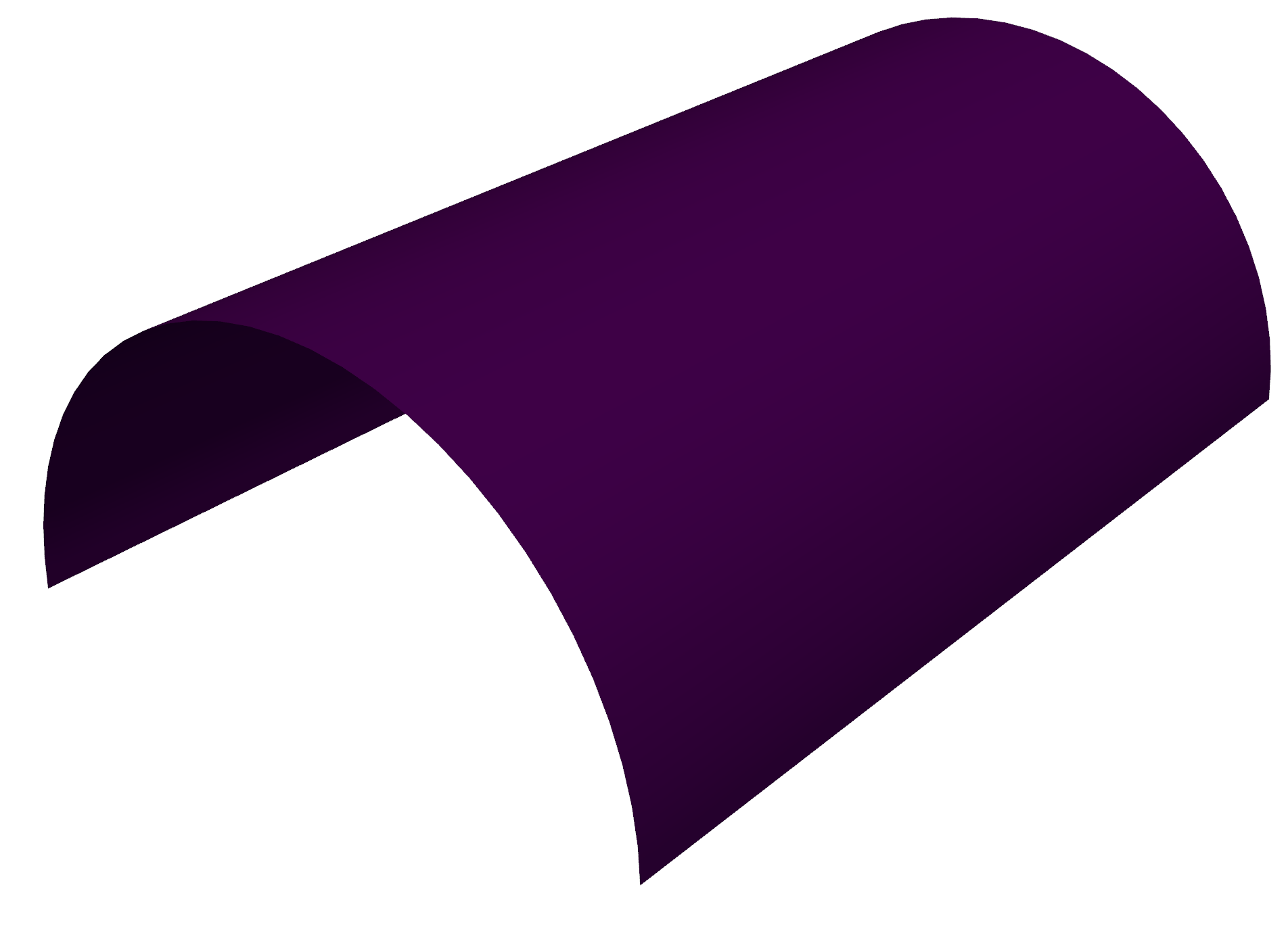}
\end{subfigure}\hspace{0.03\linewidth}
\begin{subfigure}[c]{0.31\linewidth}
  \includegraphics[width=\linewidth]{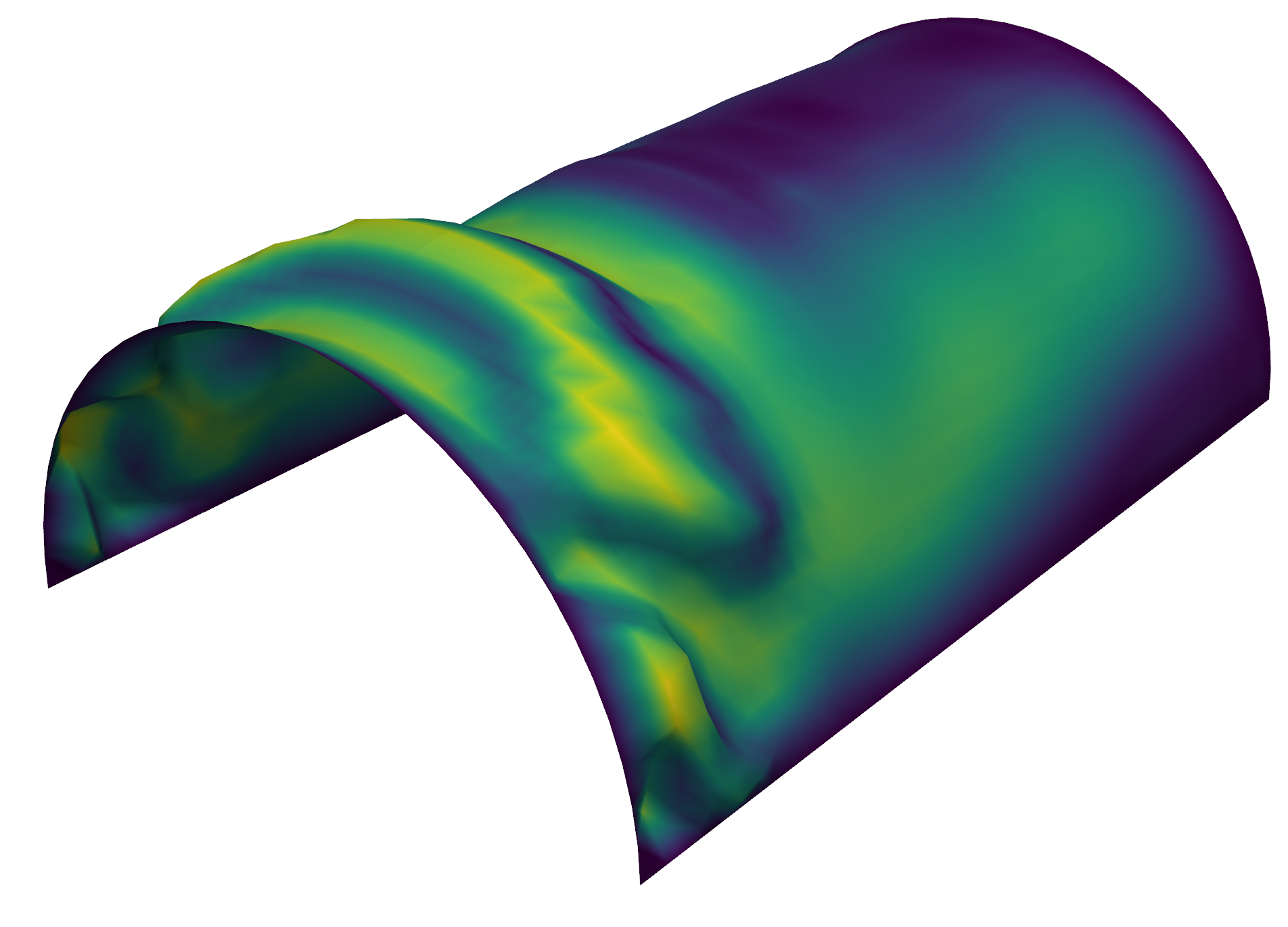}
\end{subfigure}\hspace{0.005\linewidth}
\semicylcbar{viridis}{95.2}{0,25,50,75}{Shape displacement (mm)}{3.1cm}{}

\vspace{0.7em}

\begin{subfigure}[c]{0.31\linewidth}
  \includegraphics[width=\linewidth]{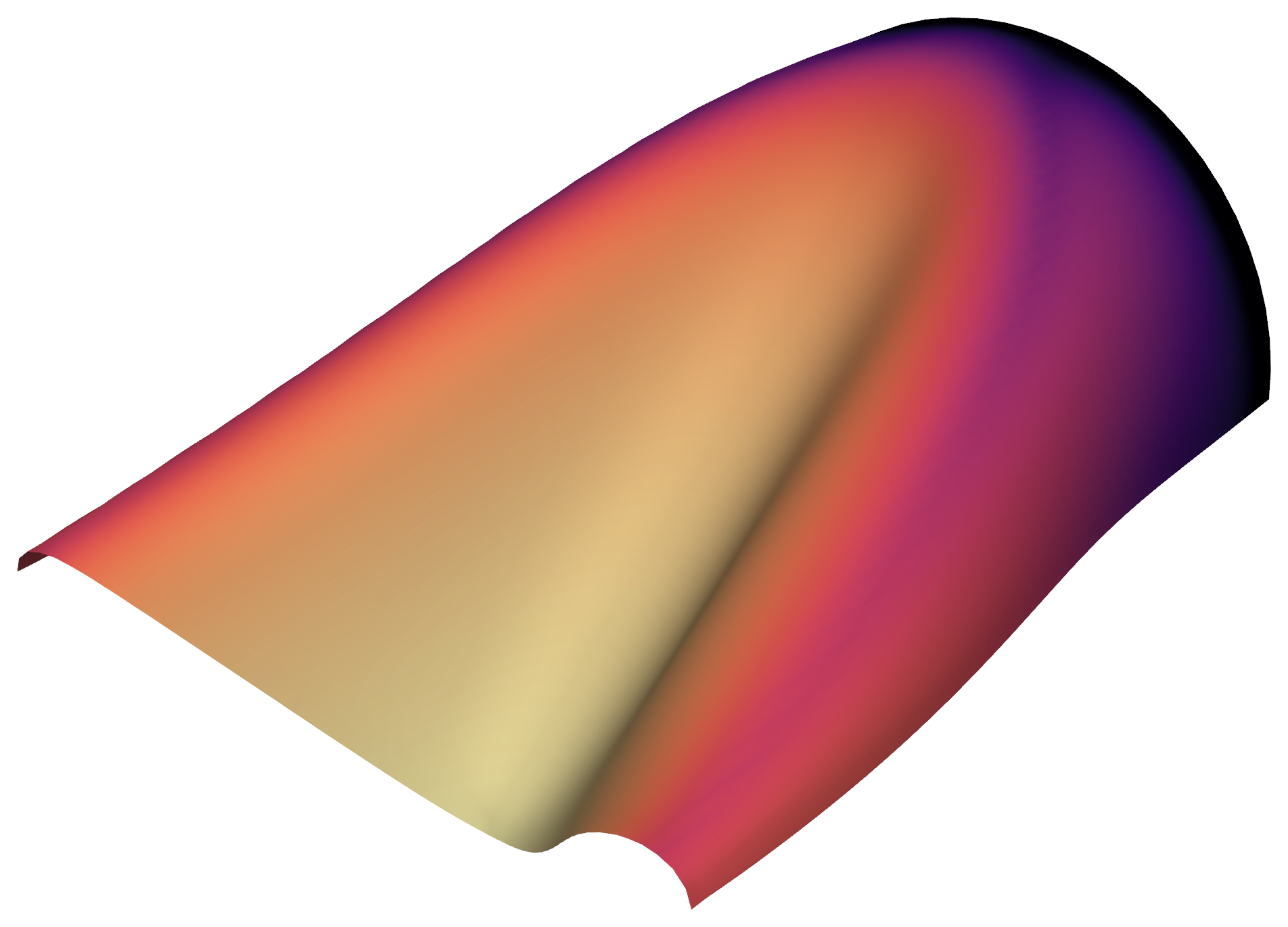}
  \caption{Initial shape}
\end{subfigure}\hspace{0.03\linewidth}
\begin{subfigure}[c]{0.31\linewidth}
  \includegraphics[width=\linewidth]{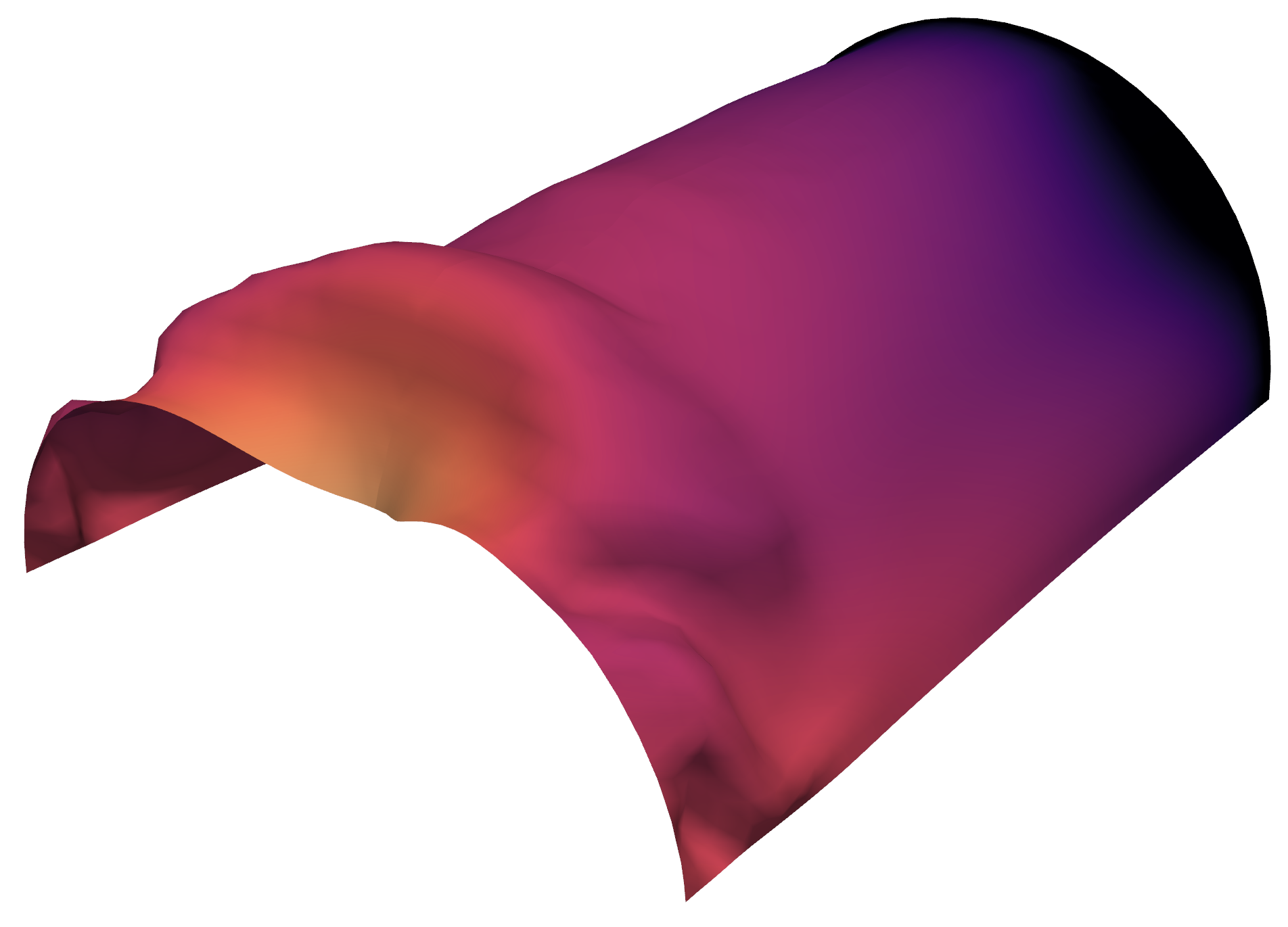}
  \caption{Optimised shape}
\end{subfigure}\hspace{0.005\linewidth}
\begin{minipage}[c]{0.15\linewidth}
\centering
\begin{tikzpicture}
\begin{axis}[
    hide axis, scale only axis, height=3.1cm, width=1pt,
    colormap name=magma, colorbar, point meta min=-3, point meta max=0.30103,
    colorbar style={
      width=0.26cm, height=3.1cm,
      ytick={-2,-1,0},
      yticklabels={$0.01$,$0.1$,$1$},
      yticklabel pos=right, ytick pos=right,
      ylabel={Displacement (m)}, ylabel near ticks,
      tick label style={font=\footnotesize}, label style={font=\footnotesize},
      axis line style={line width=0.5pt, black}, tick align=outside,
    }]
\addplot [draw=none] coordinates {(0,0) (1,1)};
\end{axis}
\end{tikzpicture}
\end{minipage}

\caption{Semi-cylinder design and loaded response, half model reflected across $x=0$.
Top row: initial and final reference surfaces coloured by shape
displacement. Bottom row: loaded configurations
coloured by displacement magnitude for (a) the initial and (b) the
100-iteration design.}
\label{fig:semicyl-opt}
\end{figure}

\paragraph{Shape optimisation}
\label{sec:results-semicyl-opt}
Shape optimising redistributes the shell mid-surface into a stiffening crease along the
crown near the loaded free end, flanked by two bulges, as visible in Figure~\ref{fig:semicyl-opt},
column (b).
The shape displacement is within the imposed \SI{\dispBudgetSemicyl}{\meter} limit
(maximum \semicylMaxDisp\,\si{\meter}). Under the \SI{2}{\kilo\newton} crown
load this reduces the squared-displacement functional by
\semicylReductionPct\,\si{\percent} relative
to the reference cylinder; equivalently the displacement $L^2$ norm falls by
\semicylAvgDeflPct\,\si{\percent}, while the surface area increases only by $\semicylAreaPct\,\si{\percent}$.
The optimisation used \semicylIters\ ROL iterations,
\semicylEvals\ objective and \semicylGrads\ gradient evaluations, and
\semicylWall\,\si{\second} of wall time on eight ranks; the replayed Newton solve
diverged at \semicylFallbacks\ of the trial shapes, each recovered by load
continuation. The run stopped at the iteration limit with gradient norm
$0.0837$.

\begin{figure}[htb!]
\centering
\begin{subfigure}[c]{0.2\linewidth}
  \includegraphics[width=\linewidth]{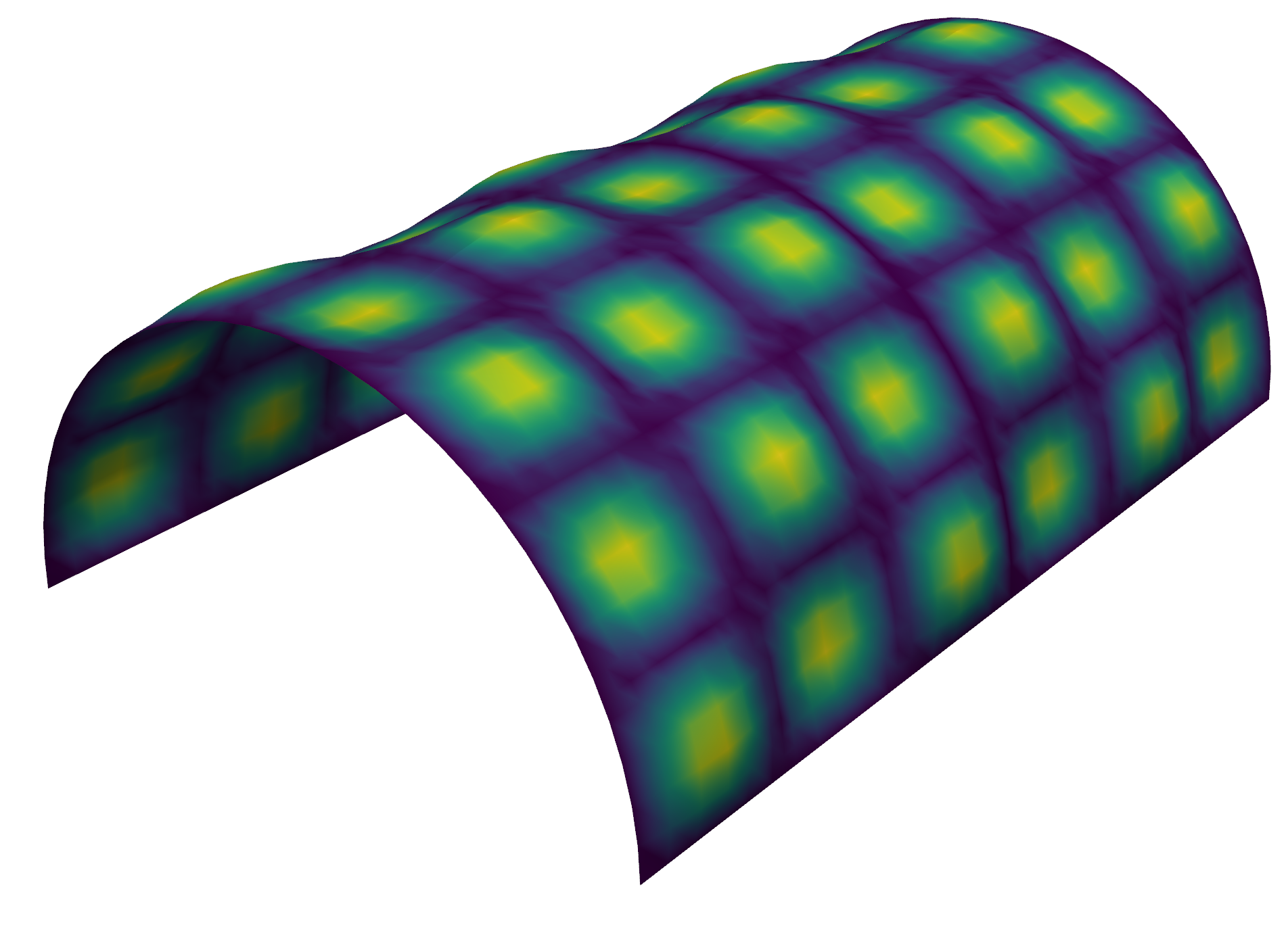}
  \caption{}
\end{subfigure}\hspace{0.002\linewidth}%
\semicylcbar{viridis}{\sinAmplitudeSemicyl}{\sinTicksSemicyl}{Perturbation (mm)}{2.6cm}{}%
\hspace{0.01\linewidth}%
\begin{subfigure}[c]{0.2\linewidth}
  \includegraphics[width=\linewidth]{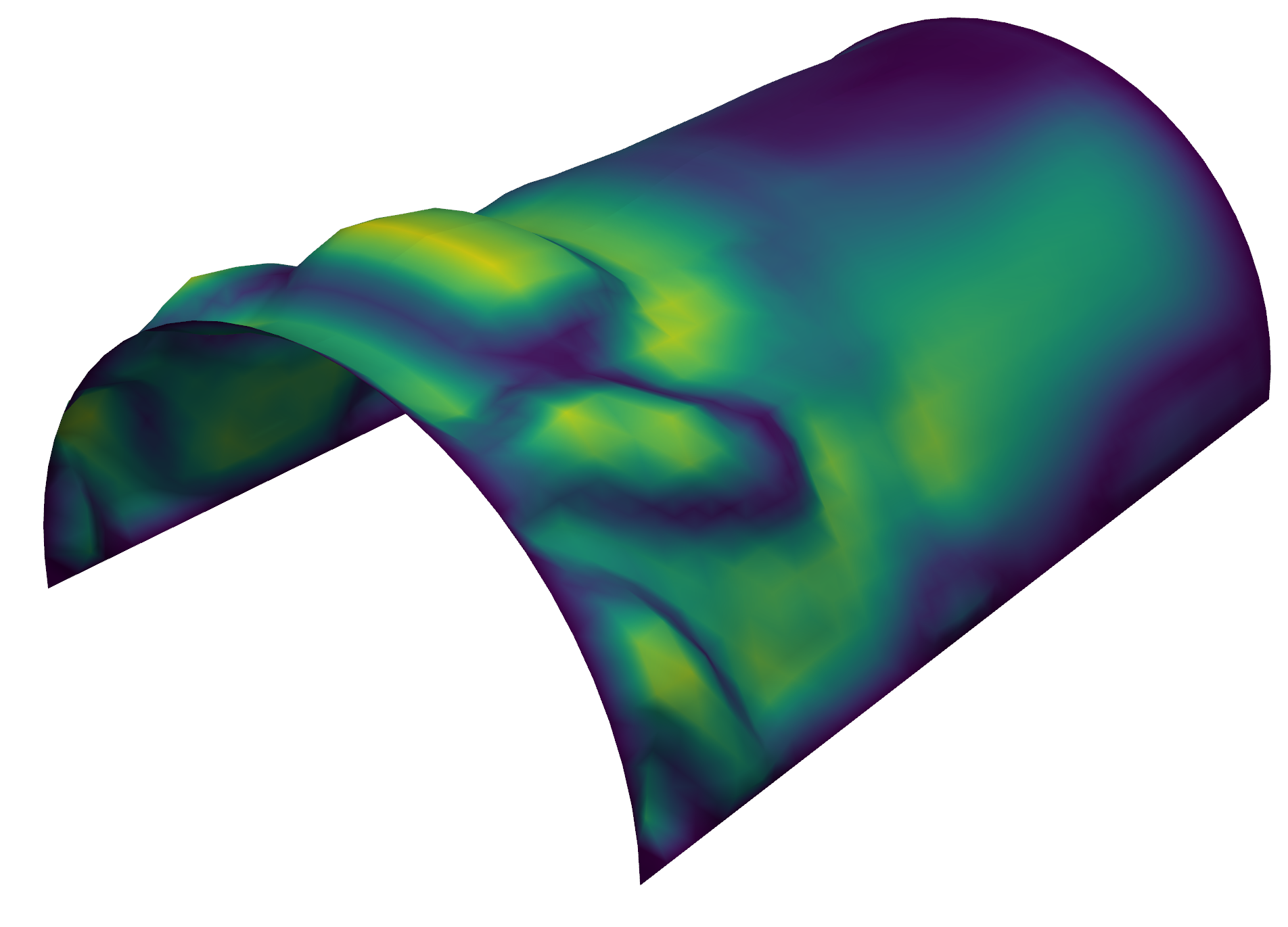}
  \caption{}
\end{subfigure}\hspace{0.01\linewidth}%
\begin{subfigure}[c]{0.2\linewidth}
  \includegraphics[width=\linewidth]{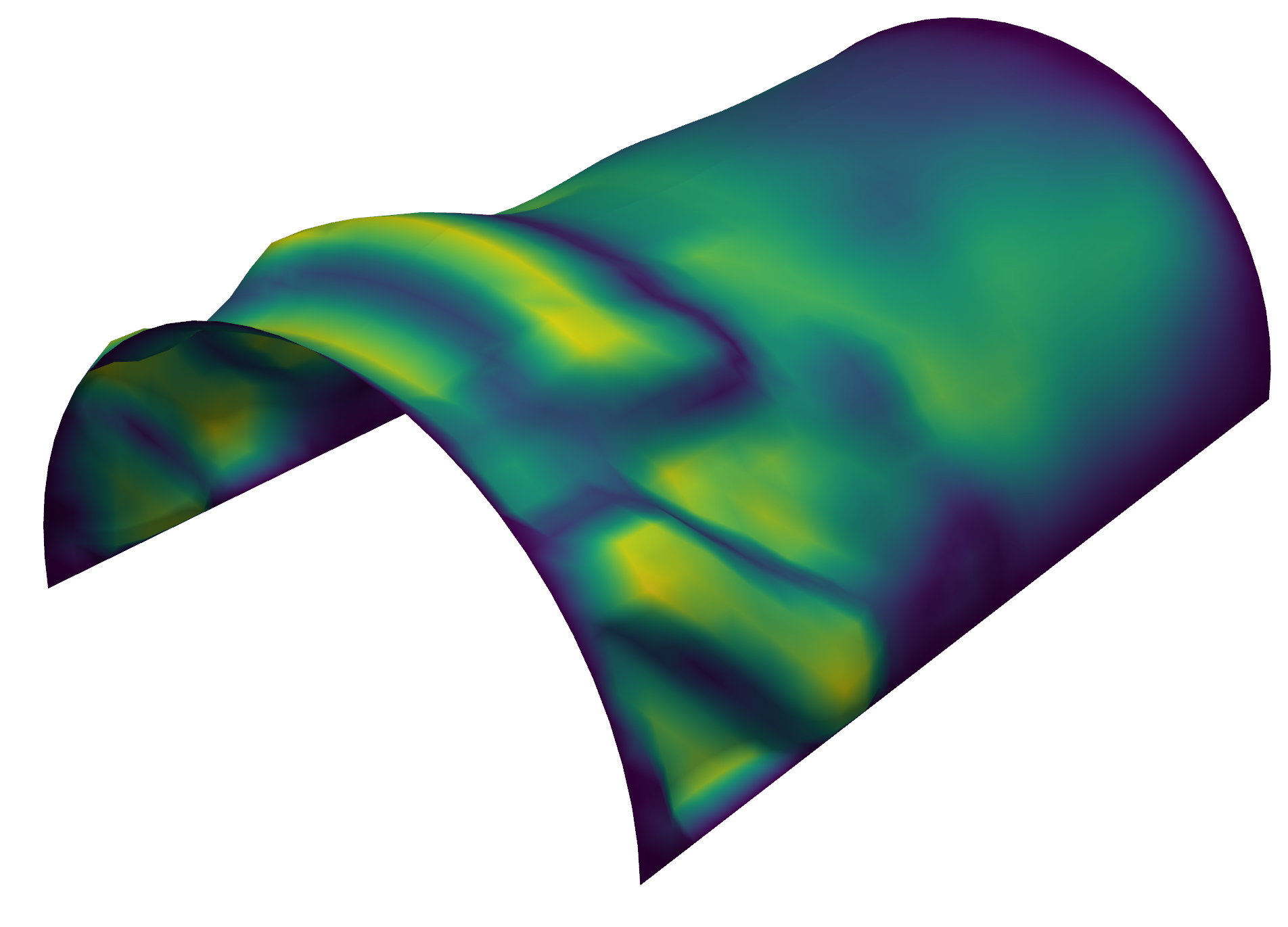}
  \caption{}
\end{subfigure}\hspace{0.002\linewidth}%
\semicylcbar{viridis}{95.2}{0,25,50,75}{Shape displacement (mm)}{2.6cm}{}
\caption{Semi-cylinder optimised from sinusoidally perturbed initial geometry: (a) the
reference cylinder perturbed by a ripple with $A=\sinStartSignSemicyl\sinAmplitudeSemicyl$\,mm, coloured by that perturbation; (b) the design optimised from
it and (c) from the mirrored case $A=\sinOppositeSignSemicyl\sinAmplitudeSemicyl$\,mm,
coloured by the shape displacement from the optimal shape obtained starting with the
unperturbed cylinder, which is displayed in the top right corner of
Figure~\ref{fig:semicyl-opt}.}
\label{fig:semicyl-starts}
\end{figure}

\paragraph{Dependence on the initial geometry and on the mesh}
\label{sec:results-semicyl-robustness}
To assess sensitivity to the initial geometry, the optimisation was
repeated with different initial guesses. Specifically, we perturbed the cylinder geometry along its
normal by adding a sinusoidal ripple of \sinModesWordsSemicyl,
($\delta\bm{X} = A\sin(\sinModesLength\pi y/L)\cos(\sinModesWidth\theta)\,\bn$ with $\theta$ measured from the
crown) that vanishes on every held edge, at \sinAmplitudeListSemicyl, that is
\sinBudgetPctSemicylSteps\,\si{\percent} of the displacement budget.
Figure~\ref{fig:semicyl-starts}(a) shows the perturbed initial shape for
$A=\sinStartSignSemicyl\sinAmplitudeSemicyl$\,mm.
The results are summarised in Table~\ref{tab:starts-semicyl}.
Several initial guesses yield similar objective values but with visibly different surfaces,
while others perform substantially worse due to an early stop of the optimisation algorithm
due to mesh distortion. (Mesh-quality promoting penalty terms could be included in the objective
to mitigate this issue.)
This
experiment evidences sensitivity to starting geometry.
We also performed a refinement experiment. For the $A=\pm20$\,mm case,
the state mesh was uniformly refined from $4{,}096$ to $16{,}384$ triangles,
while the control mesh remained unchanged ($1{,}024$ triangles).
Working with finer state meshes leads to only minor changes with respect to their
coarse-mesh counterparts, as evidenced in Table~\ref{tab:refinement}.

\begin{table}[htb!]\centering\small
\caption{Semi-cylinder: initial-geometry sweep at the final recorded iteration. Objective eductions and initial and final area changes
$\Delta A_0$ and $\Delta A_f$ refer to the unperturbed cylinder. Root-mean squared (RMS) and maximum (Max) distances
compare the final mid-surface with the unperturbed-start final design and
quantify variation across optimised designs.
Maximum nodal shape displacement from the initial (i.e., unoptimised) cylinder $\max|q_f|$ indicates whether the prescribed displacement budget of \SI{100}{\milli\meter} is satisfied.}
\label{tab:starts-semicyl}
\begin{tabular}{@{}rrrrrrrr@{}}\toprule
$A$ & Iter. & Reduction & $\Delta A_0$ & $\Delta A_f$ & $\max|q_f|$ & RMS & Max.\\
(mm) & & (\%) & (\%) & (\%) & (mm) & (mm) & (mm)\\\midrule
$+2$ & 54 & 72.86 & 0.00 & 0.58 & 93.0 & 52.2 & 168.1 \\
$-2$ & 100 & 98.46 & 0.00 & 1.74 & 94.4 & 17.2 & 71.3 \\
$+5$ & 39 & 65.31 & 0.03 & 0.37 & 94.2 & 51.8 & 171.1 \\
$-5$ & 100 & 98.60 & 0.03 & 1.76 & 91.3 & 22.2 & 82.7 \\
$+10$ & 100 & 98.50 & 0.10 & 1.84 & 94.8 & 34.1 & 143.4 \\
$-10$ & 61 & 64.71 & 0.10 & 0.26 & 96.4 & 51.2 & 169.9 \\
$+20$ & 100 & 98.60 & 0.41 & 1.81 & 95.4 & 33.9 & 149.5 \\
$-20$ & 100 & 98.62 & 0.41 & 1.79 & 92.8 & 20.9 & 83.8 \\
$+40$ & 100 & 98.51 & 1.64 & 1.83 & 95.2 & 36.5 & 149.1 \\
$-40$ & 100 & 98.75 & 1.64 & 1.80 & 91.0 & 17.8 & 79.9 \\
\bottomrule\end{tabular}\end{table}

\begin{table}[htb!]\centering\small
\caption{Semi-cylinder: state refinement experiment optimising from perturbed design. The ``coarse'' state mesh has $4{,}096$ triangles while
the ``fine'' state mesh as $16{,}384$. In both cases, the control mesh has $1{,}024$ triangles.
Each fine design is compared directly with its matching coarse design. The differences in optimised designs obtained with ``coarse''
and ``fine'' state meshes are minor.}
\label{tab:refinement}
\begin{tabular}{@{}rrrrrr@{}}\toprule
$A$ (mm) & Iter. (c/f) & Coarse (\%) & Fine (\%) & RMS (mm) & Max. (mm)\\\midrule
$+20$ & 100/100 & 98.60 & 98.64 & 3.11 & 16.65 \\
$-20$ & 100/100 & 98.62 & 98.62 & 4.68 & 25.44 \\
\bottomrule\end{tabular}\end{table}

\subsection{Sheet-metal bracket: the COMSOL benchmark}
\label{sec:results-bracket}

\paragraph{Problem definition}
The COMSOL ``Shape Optimization of a Shell'' benchmark \cite{comsol_shell_shapeopt} is an
L-shaped steel shell with a filleted corner, loaded by a total force of \SI{10}{\kilo\newton}
in the $y$-direction distributed uniformly over the free edge and clamped on the
opposite edge; see Figure~\ref{fig:bracket-shape-sketch}.
The bracket is optimised for maximum stiffness, which is equivalent to minimising the total
elastic strain energy $J_{\mathrm{br}}$ described in Section~\ref{sec:shapeopt}.
Geometry and material values are collected in Table~\ref{tab:bracket-setup}.
As in the benchmark, shape updates can only slide the clamped edge within its own plane
(a roller, $q_y=0$) while the remaining external edges are held fixed
($\bq=\bm{0}$). The admissible shape change is limited to
$\delta_{\max}=\SI{0.05}{\meter}$, which is implemented using the soft penalty with smoothing
parameter $\varepsilon=\SI{2.5}{\milli\meter}$. The second budget is the admissible relative
increase of surface area, $a_{\max}=+\areaBudgetBracket\,\si{\percent}$, which is the value
attained in the benchmark design. The problem is symmetric about the
half-height mid-plane $z=\tfrac12 H$, so the half model is optimised with a symmetry-plane
condition and mirrored. The control mesh has $6{,}080$ triangles (the $380$-triangle
half mesh refined twice) and the state mesh $24{,}320$ triangles
($319{,}605$ degrees of freedom); the Helmholtz length scale is
$\ell=\SI{0.6}{\meter}$ and the load continuation uses five increments.

\begin{table}[htb!]
\centering
\small
\caption{Sheet-metal bracket (COMSOL ``Shape Optimization of a Shell''
\cite{comsol_shell_shapeopt}) benchmark: geometry, material, and conditions.}
\label{tab:bracket-setup}

\begin{tabular}{@{}lll@{}}
\toprule
\multicolumn{3}{@{}l}{\textbf{Sheet-metal bracket}}\\
\midrule

\textit{Geometry}
&
\textit{Material}
&
\textit{Conditions}
\\[2pt]

\begin{tabular}[t]{@{}ll@{}}
Footprint $L$
  & \SI{1.0}{\meter} \\
Height $H$
  & \SI{0.5}{\meter} \\
Fillet radius $R$
  & \SI{0.3}{\meter} \\
Thickness $t$
  & \SI{10}{\milli\meter}
\end{tabular}
&
\begin{tabular}[t]{@{}ll@{}}
Young's modulus $E$
  & \SI{200}{\giga\pascal} \\
Poisson ratio $\nu$
  & \num{0.3}
\end{tabular}
&
\begin{tabular}[t]{@{}ll@{}}
Applied load $P$
  & \SI{10}{\kilo\newton} \\
Shape budget $\delta_{\max}$
  & \SI{50}{\milli\meter} \\
Area budget $a_{\max}$
  & $+\areaBudgetBracket\,\si{\percent}$\\
Helmholtz scale $\ell$
  & \SI{0.6}{\meter} \\
Control/state mesh
  & $6{,}080/24{,}320$ triangles
\end{tabular}
\\

\bottomrule
\end{tabular}
\end{table}

\begin{figure}[htb!!]
\centering
\begin{subfigure}[c]{\linewidth}
\centering
\input{bracket_setup_sketch.tex}
\caption{}
\end{subfigure}
\caption{Reference geometry, boundary conditions, loads, and fixed/clamped (solid/dotted blue)
optimisation boundaries for the bracket test case.}
\label{fig:bracket-shape-sketch}
\end{figure}

\begin{figure}[htb!!]
\centering
\begin{subfigure}[c]{0.22\linewidth}
  \includegraphics[width=\linewidth]{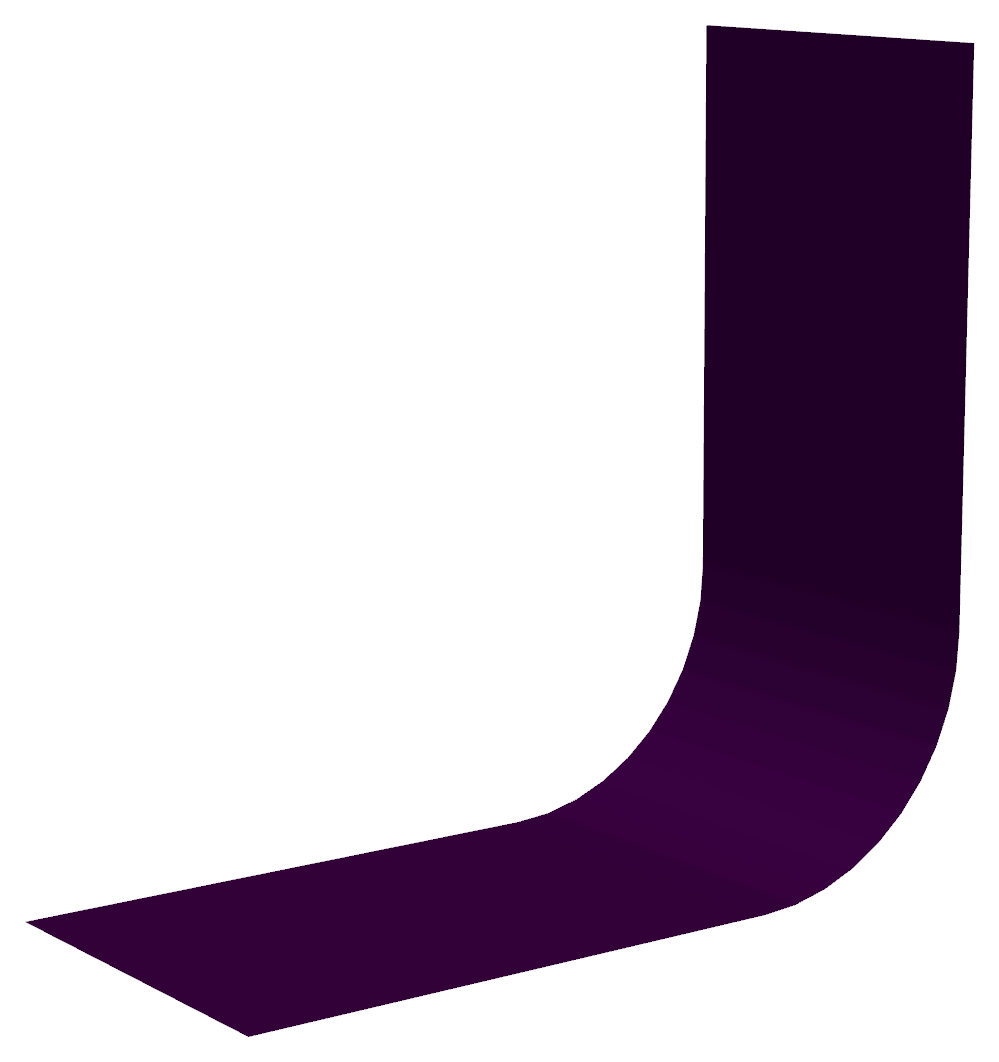}
\end{subfigure}\hspace{0.03\linewidth}
\begin{subfigure}[c]{0.22\linewidth}
  \includegraphics[width=\linewidth]{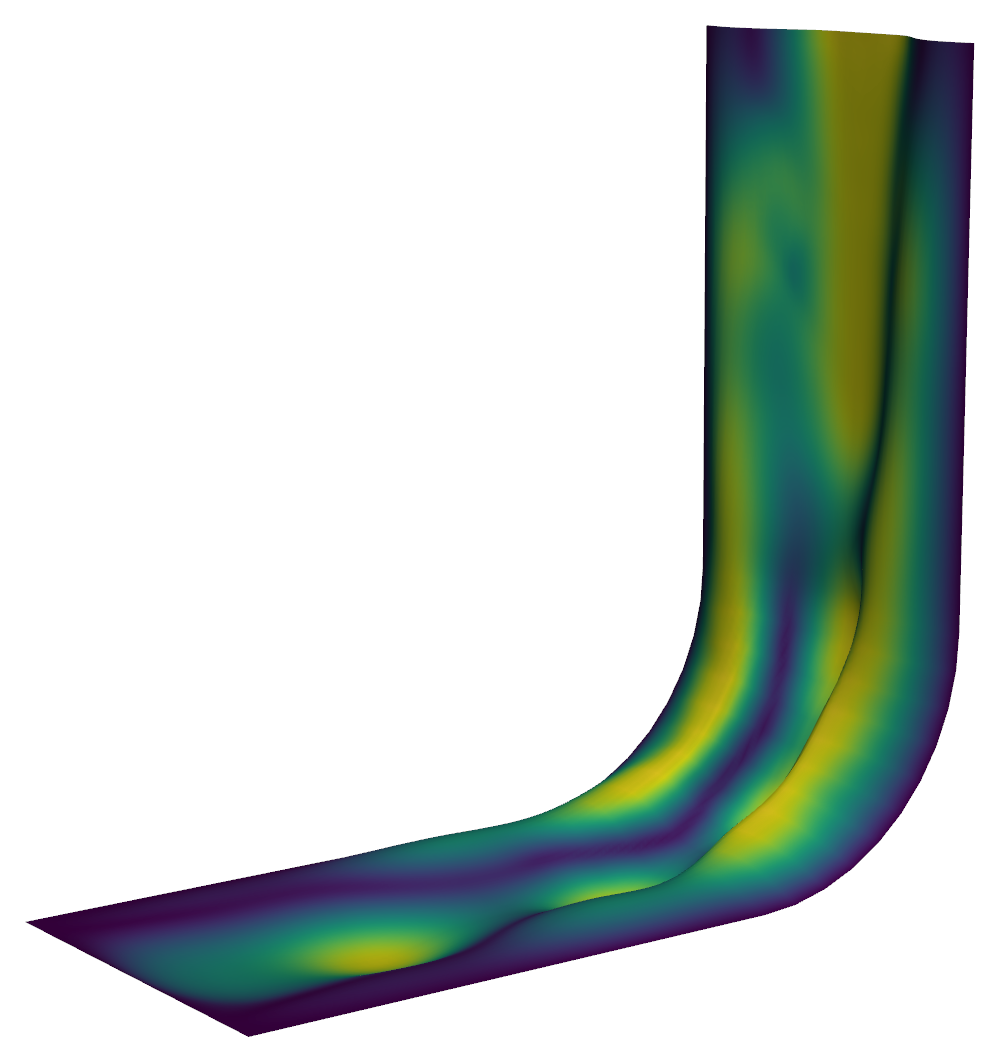}
\end{subfigure}\hspace{0.005\linewidth}
\semicylcbar{viridis}{49.9}{0,10,20,30,40}{Shape displacement (mm)}{3.1cm}{}

\vspace{0.7em}

\begin{subfigure}[c]{0.22\linewidth}
  \includegraphics[width=\linewidth]{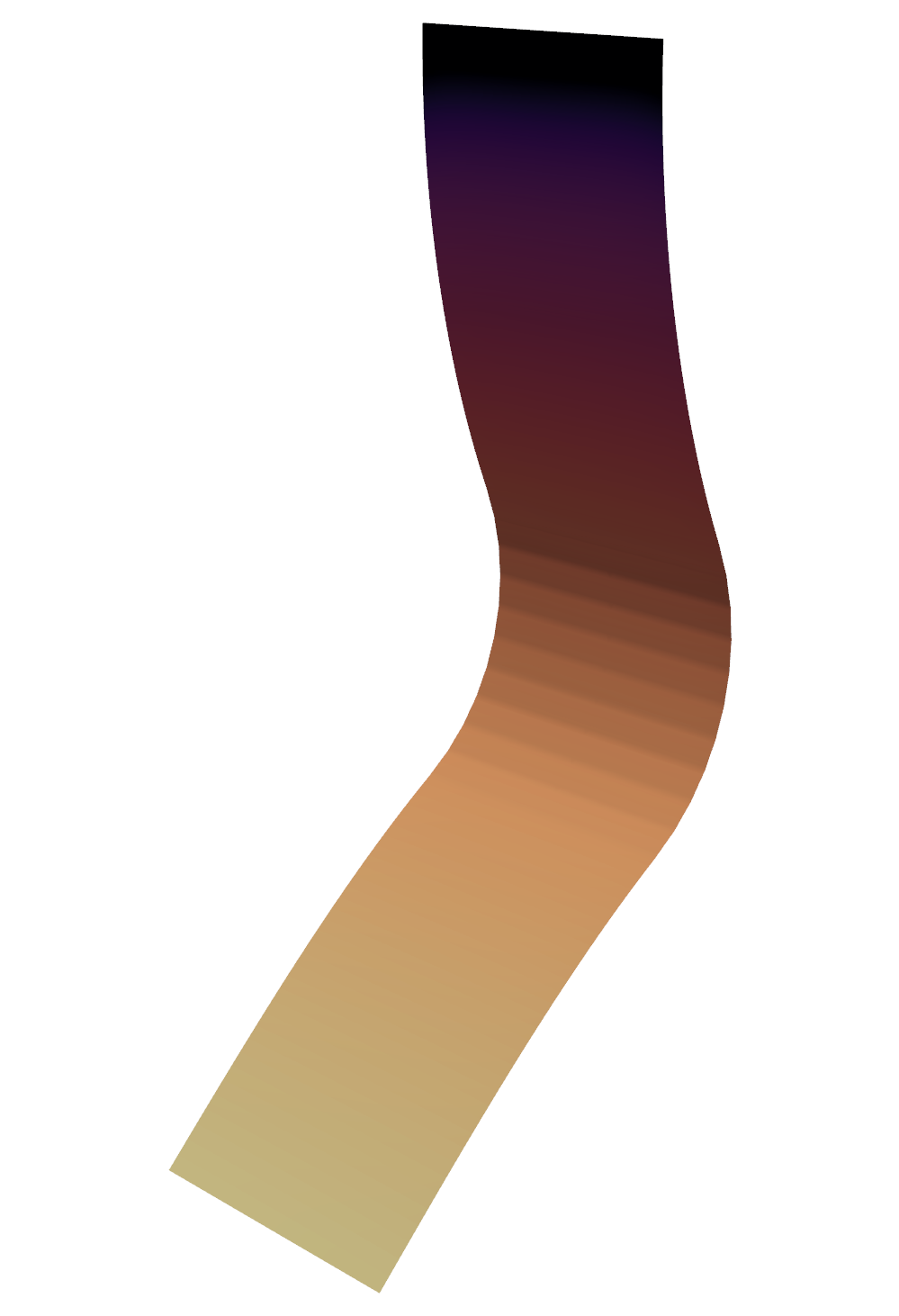}
  \caption{Initial shape}
\end{subfigure}\hspace{0.03\linewidth}
\begin{subfigure}[c]{0.22\linewidth}
  \includegraphics[width=\linewidth]{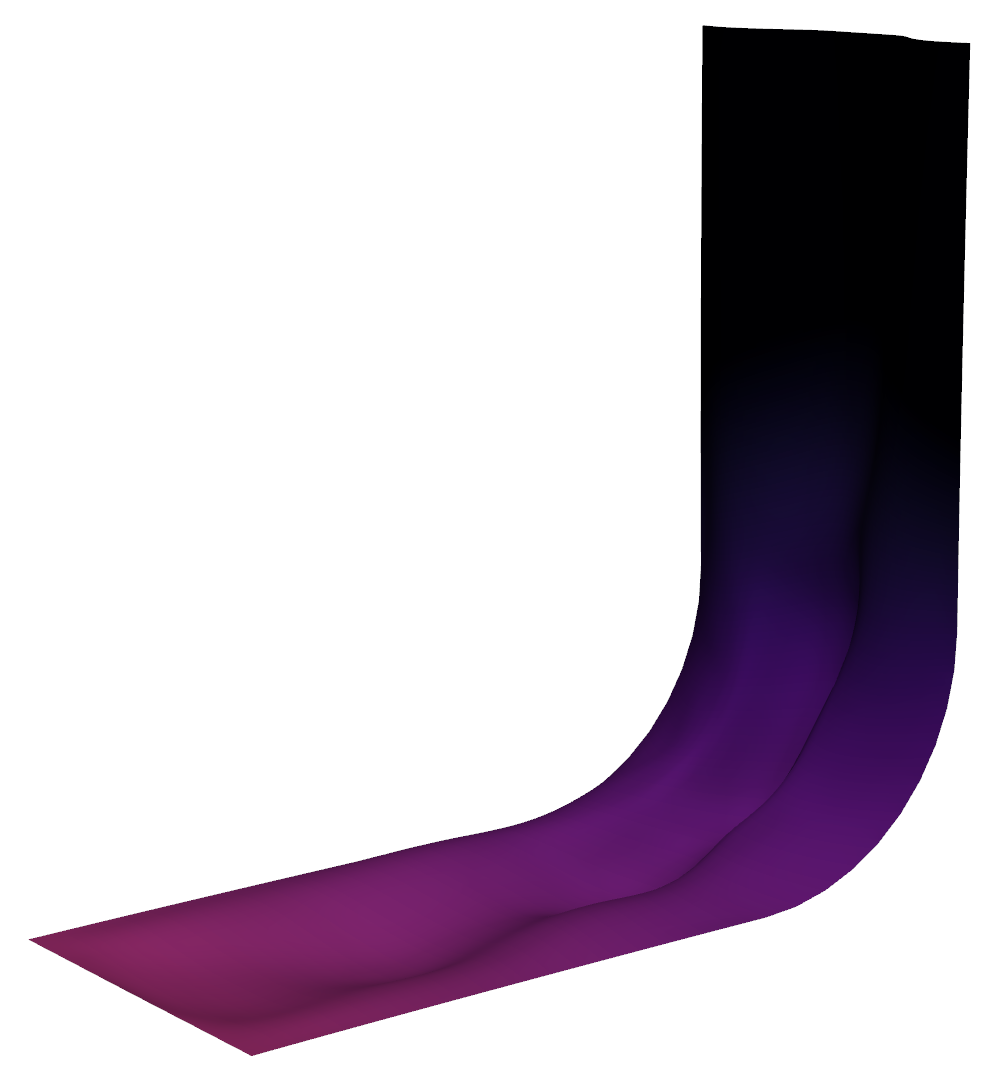}
  \caption{Optimised shape}
\end{subfigure}\hspace{0.005\linewidth}
\begin{minipage}[c]{0.15\linewidth}
\centering
\begin{tikzpicture}
\begin{axis}[
    hide axis, scale only axis, height=3.1cm, width=1pt,
    colormap name=magma, colorbar, point meta min=-3, point meta max=0,
    colorbar style={
      width=0.26cm, height=3.1cm,
      ytick={-2,-1,0},
      yticklabels={$0.01$,$0.1$,$1$},
      yticklabel pos=right, ytick pos=right,
      ylabel={Displacement (m)}, ylabel near ticks,
      tick label style={font=\footnotesize}, label style={font=\footnotesize},
      axis line style={line width=0.5pt, black}, tick align=outside,
    }]
\addplot [draw=none] coordinates {(0,0) (1,1)};
\end{axis}
\end{tikzpicture}
\end{minipage}
\caption{Bracket design and loaded response, half model reflected across its mid-plane.
Top row: initial and final reference surfaces coloured by shape displacement.
Bottom row: loaded configurations coloured by displacement magnitude for
(a) the initial and (b) the 100-iteration design.}
\label{fig:bracket-shape}
\end{figure}

\paragraph{Shape optimisation}
The optimiser moves material away from the mid-plane, forming two stiffening
corrugations along the legs and around the fillet (Figure~\ref{fig:bracket-shape})
that recruit bending stiffness from the initially flat regions. The design update reduces the
elastic strain energy by \bracketReductionPct\,\si{\percent}
and essentially makes full use (but does not exceed) the \SI{0.05}{\meter} displacement and
$+\areaBudgetBracket\,\si{\percent}$ area budgets.
For comparison, COMSOL achieved \SI{89}{\percent} reduction with a
\SI{0.05}{\meter} displacement bound. Surface area was not constrained in that benchmark, and
the reported \SI{9}{\percent} increase is just measured after optimisation.
We impose an additional soft area budget to make the attained areas comparable.
The reduction is \bracketGainOverComsol\ percentage points greater, but the
shell formulation, design space and constraint treatment differ, so these are
benchmark outcomes rather than identical optimisation problems. The run stopped
at the iteration limit (\bracketIters) with gradient norm $0.174$. The optimisation
algorithm used \bracketEvals\ objective and
\bracketGrads\ gradient evaluations and \bracketWall\,\si{\second} on eight
ranks, with \bracketFallbacks\ replays recovered by load continuation. The
final state mesh has minimum triangle quality \meshQualityBracket.

Dropping the area budget allows the optimiser to achieve even higher reduction.
For example, we can obtain
\unboundedBracketReductionPct\,\si{\percent} reduction at iteration
\unboundedBracketReportedIter\ with $\unboundedBracketAreaPct\,\si{\percent}$ of
added surface (not displayed). Further reduction is achievable by adding ridges of ever shorter wavelength
(while staying within the displacement budget; not displayed). For example, the algorithm
achieves \unboundedBracketFinalReductionPct\,\si{\percent} reduction with
$\unboundedBracketFinalAreaPct\,\si{\percent}$ of added area at the iteration limit, where the
state mesh reaches the minimum triangle quality allowed.

\begin{figure}[htb!]
\centering
\begin{subfigure}[c]{0.2\linewidth}
  \includegraphics[width=\linewidth]{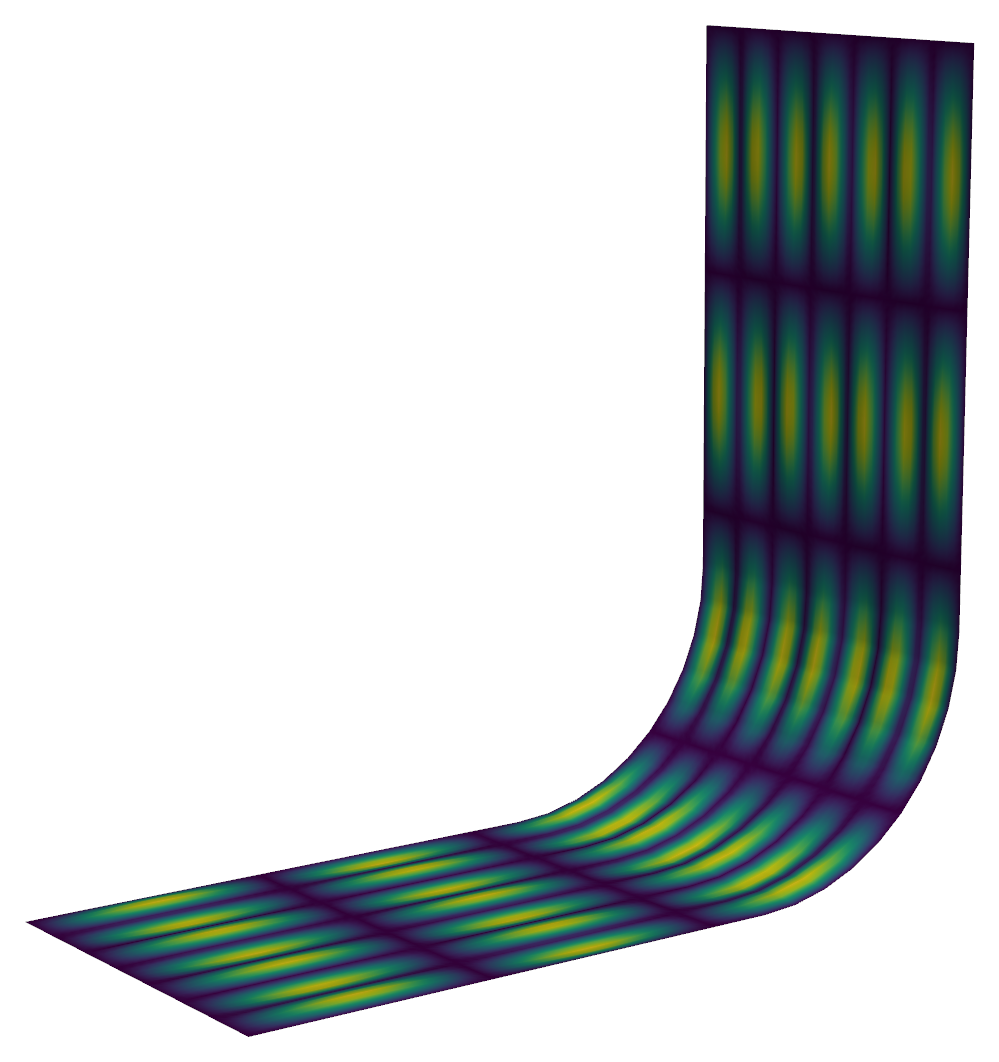}
  \caption{}
\end{subfigure}\hspace{0.002\linewidth}%
\semicylcbar{viridis}{\sinAmplitudeBracket}{\sinTicksBracket}{Perturbation (mm)}{2.6cm}{}%
\hspace{0.01\linewidth}%
\begin{subfigure}[c]{0.2\linewidth}
  \includegraphics[width=\linewidth]{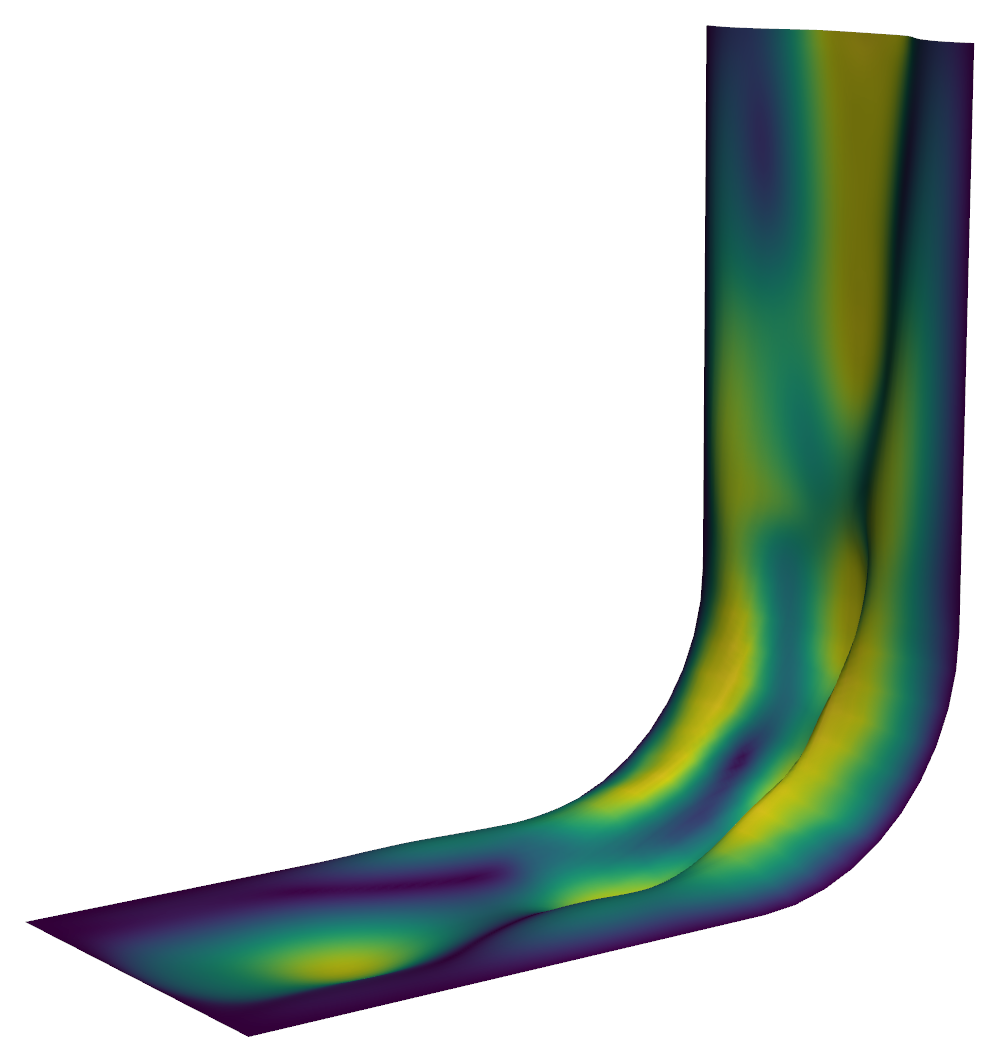}
  \caption{}
\end{subfigure}\hspace{0.01\linewidth}%
\begin{subfigure}[c]{0.2\linewidth}
  \includegraphics[width=\linewidth]{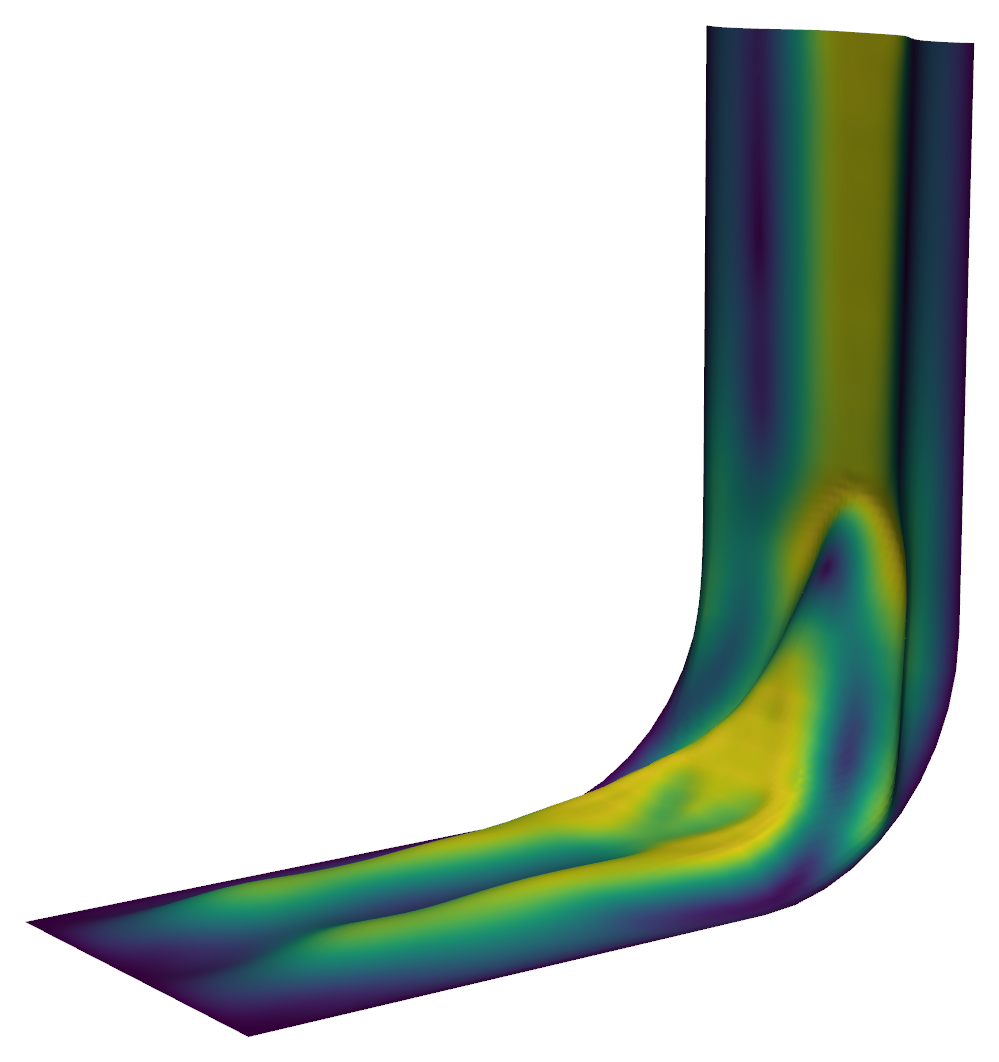}
  \caption{}
\end{subfigure}\hspace{0.002\linewidth}%
\semicylcbar{viridis}{49.9}{0,10,20,30,40}{Shape displacement (mm)}{2.6cm}{}
\caption{Bracket optimised from sinusoidally perturbed initial geometry:
(a) the flat bracket perturbed by a ripple with
$A=\sinStartSignBracket\sinAmplitudeBracket$\,mm, coloured by
that perturbation; (b) the design optimised from it and (c) from the mirrored
case $A=\sinOppositeSignBracket\sinAmplitudeBracket$\,mm, coloured by the shape
displacement from
the optimal shape obtained starting with the unperturbed bracket,
which is displayed in the top right corner of Figure~\ref{fig:bracket-shape}.}
\label{fig:bracket-starts}
\end{figure}

\paragraph{Dependence on the initial geometry}
Similarly to the previous test case, we assess sensitivity to the initial geometry
by repeating the optimisation with different initial guesses.
The flat bracket was perturbed along its normal by a sinusoidal ripple of
\sinModesWordsBracket
($\delta\bm{X} = A\sin(\sinModesLength\pi s/S)\sin(\sinModesWidth\pi z/W)\,\bn$ with $s$ the arc length from the clamped edge and $S$ the total arc length) that vanishes on every held edge, at
\sinAmplitudeListBracket, that is \sinBudgetPctBracketSteps\,\si{\percent}
of the displacement budget. Figure~\ref{fig:bracket-starts}(a) shows the perturbed
initial shape for $A = 5$\,mm. The results are summarised in Table~\ref{tab:starts-bracket}.
All initial guesses yield similar objective values, albeit most designs present
visible and quantifiable differences.
The $+20$\,mm start has an area increase of $9.12\%$, which is slightly above the
nominal $9\%$ soft budget. This is possible because penalisation discourages but does
not exclude constraint violation.
These results evidence sensitivity to the initial guess.

\begin{table}[htb!]\centering\small
\caption{Bracket: initial-geometry sweep at the final recorded iteration. Objective reductions and 
initial and final area changes $\Delta A_0$ and $\Delta A_f$ refer to the unperturbed L-shape design.
Root-mean squared (RMS) and maximum (Max) distances compare the final mid-surface with the
unperturbed-start final design. Maximum nodal shape displacement from the initial (i.e., unoptimised) cylinder $\max|q_f|$ indicates whether the prescribed displacement budget of \SI{100}{\milli\meter} is satisfied.}
\label{tab:starts-bracket}
\begin{tabular}{@{}rrrrrrrr@{}}\toprule
$A$ & Iter. & Reduction & $\Delta A_0$ & $\Delta A_f$ & $\max|q_f|$ & RMS & Max.\\
(mm) & & (\%) & (\%) & (\%) & (mm) & (mm) & (mm)\\\midrule
$+5$ & 100 & 94.59 & 0.62 & 8.94 & 50.6 & 3.0 & 13.6 \\
$-5$ & 100 & 93.88 & 0.62 & 8.92 & 50.8 & 31.0 & 95.4 \\
$+10$ & 100 & 94.44 & 2.42 & 8.92 & 50.8 & 27.1 & 97.0 \\
$-10$ & 100 & 93.25 & 2.42 & 8.94 & 51.4 & 18.1 & 77.6 \\
$+20$ & 100 & 93.47 & 9.12 & 8.93 & 54.3 & 24.0 & 87.5 \\
$-20$ & 100 & 93.19 & 9.12 & 8.94 & 51.3 & 18.4 & 87.2 \\
\bottomrule\end{tabular}\end{table}

\subsection{Impact of nonlinear modelling}
\label{sec:results-linear}
For both test cases, we also performed optimisation with
linearised kinematics constraints
(first-order strains in the displacement and the director increment).
The resulting designs were then assessed with the geometrically nonlinear model
at the design load.
These comparisons show that the linear model can substantially mispredict the response.
In the semi-cylinder case, the predicted squared-displacement
functional of the initial design at \SI{2}{\kilo\newton} is
\linRefFactorSemicyl\ times smaller than the nonlinear value (\SI{0.098}{\metre\tothe{4}}
against \SI{1.594}{\metre\tothe{4}}, the nonlinear shell being far softer at
this load).
In the bracket test case, the strain energy of the initial design is
\linRefRatioBracket\ times larger than the nonlinear value (\SI{3395}{\joule} against
\SI{1065}{\joule}).

With optimisation, the linear model reduces the objective values by
\linOwnReductionPctSemicyl\si{\percent} and \linOwnReductionPctBracket\si{\percent}
(using a maximum geometry perturbation of \linMaxDispSemicyl\,\si{\meter}
and \linMaxDispBracket\,\si{\meter}, respectively).
Following this significant reduction, the displacements become sufficiently small to justify a
linear modelling assumption. Indeed, when assessed against the nonlinear model, the
differences are small. Specifically, the linear and nonlinear models predict objective values of
\SI{0.005}{\metre\tothe{4}} and \SI{0.012}{\metre\tothe{4}} for the optimised semi-cylinder,
and \SI{\linFinalRawBracket}{\joule} and \SI{\linNonlinearRawBracket}{\joule}
for the optimised bracket, respectively.
This indicates that the considerably cheaper linear model can still identify effective stiffening mechanisms,
despite providing inaccurate quantitative predictions of the mechanical response during the initial optimisation
iterations, when deformations exceed the limits of the small-deformation assumption.

\subsection{Verification and cost}
\label{sec:results-verification}

\emph{Derivative validation.} A Taylor test of the shell objective $J$ (without
penalty terms) at the reference shape $\bq=0$ performed along a smooth random
control direction $\bdq$ (with a maximum nodal displacement of \SI{0.05}{\meter})
shows second-order convergence of the linear approximation $s\mapsto
J(\bq)+sdJ(\bdq)$ of $J(\bq+s\bdq)$, with observed minimum rates of
\taylorRateSemicyl\ (semi-cylinder, \SI{2}{\kilo\newton}) and
\taylorRateBracket\ (bracket).  This validates the correctness of the shape
derivative $dJ$, which is obtained via algorithmic differentiation.  Adding the
Hessian action to form the quadratic model $s\mapsto J(\bq)+sdJ(\bdq)+(s^2/2)
d^2J(\bdq, \bdq)$ increases the convergence rates to \taylorRateSecondSemicyl\
and \taylorRateSecondBracket, which corresponds to the expected third-order convergence.
This validates the correctness of the Hessian action. Penalty terms were tested
separately and also achieved the expected convergence rates, which validates
the correctness of their algorithmic differentiation.

\emph{Mesh quality.} Complementing the surface-Jacobian feasibility check of
Section~\ref{sec:shapeopt}, the minimum normalised triangle quality
$m_q=4\sqrt{3}A/(a^2+b^2+c^2)$ of the state mesh was tracked along every
baseline optimisation. At the reported final iterate the semi-cylinder mesh has
$m_q=\meshQualitySemicyl$ and a minimum angle of \meshAngleSemicyl$^\circ$,
while the bracket mesh has $m_q=\meshQualityBracket$ and a minimum angle of
\meshAngleBracket$^\circ$ (compared to $m_q=0.86$ and minimum angle of
$43^\circ$ on the initial meshes). The mesh quality remains acceptable,
although it deteriorated as the crease and the corrugations stretch
the facets. We note that the area budget also helps preserve mesh
quality by discouraging the formation of additional folds.

\emph{Optimiser.} As explained in \ref{sec:shapeopt}, the experiments reported thus far
were obtained with a trust-region method with a Steihaug--Toint truncated
conjugate-gradient solver and L-BFGS quadratic model.
Fireshape also implements shape-Hessian-vector actions, which enable a more
accurate second-order model at the increased computational cost of an additional
tangent-linear state solve and a second-order adjoint solve per Hessian-vector evaluation.
When applied to the bracket test case,  the second-order-adjoint Hessian
achieves a \hessianReductionPct\,\si{\percent} objective reduction at iteration
\hessianReportedIter, requiring \hessianGrads\ gradient evaluations and
\hessianActions\ accumulated inner CG iterations. The final mesh quality is \hessianQuality.
This illustrates the availability of alternative second-order strategies within the same framework,
although these results alone do not establish their relative computational efficiency.

\emph{Computational cost.} Table~\ref{tab:cost} reports the cost of one design
evaluation cycle of the shell objective at the production resolution. Medians and
full ranges are computed over nine cycles (2--10) of shape perturbation, objective
evaluation and gradient evaluation, excluding cycles 0 and 1.
The continuation means and optimisation
times come from the area-bounded baselines.
Optimisation-loop time includes
replays, derivative evaluations, fallbacks and output, but excludes initial
setup/continuation and tape creation (about $112+21$\,s for the semi-cylinder
and $175+18$\,s for the bracket).

\emph{Software.} The baseline summaries record Firedrake
\firedrakeVersion, Fireshape commit \fireshapeCommitMain\ (September 2026)
and ROL through pyroltrilinos 0.5.6. The supplementary manifest records the
available per-run provenance and the local image identifier inspected for the
additional verification. Older timing runs have their own recorded provenance.
An unavailable historical dependency commit is left unspecified rather than
inferred from the current installation.

\begin{table}[htb!]
\centering
\caption{Shell evaluation costs at production resolution on eight MPI
ranks (64-core AMD Ryzen Threadripper PRO 9985WX): geometry recovery, objective
by tape replay (warm-started Newton solve), adjoint shape derivative and
Helmholtz Riesz map: {median [minimum, maximum] over the nine cycles 2--10.} Continuation
is the mean per fallback in the baseline; the final row is optimisation-loop time.
Counts refer to the complete baseline runs. Forward-solve attempts include
initialisation, unsuccessful solves and every load-continuation step; each
adjoint evaluation computes one shell-objective gradient.}
\label{tab:cost}
\begin{tabular}{@{}lrr@{}}
\toprule
 & Semi-cylinder & Bracket\\
 & ($54{,}213$ DOFs) & ($319{,}605$ DOFs)\\
\midrule
Geometry recovery (\si{\second})              & {\costGeomSemicyl} & {\costGeomBracket}\\
Shell objective evaluation (\si{\second}) & {\costFwdSemicyl} & {\costFwdBracket}\\
Adjoint shape derivative (\si{\second})       & {\costAdjSemicyl} & {\costAdjBracket}\\
Riesz map (\si{\second})                       & {\costRieszSemicyl} & {\costRieszBracket}\\
\midrule
Load continuation, when needed (\si{\second}) & \costContSemicyl & \costContBracket\\
Optimisation iterations                      & \semicylIters & \bracketIters\\
Forward nonlinear solve attempts             & 853 & 142\\
Adjoint evaluations                          & \semicylGrads & \bracketGrads\\
Optimisation loop (\si{\second})             & \semicylWall & \bracketWall\\
\bottomrule
\end{tabular}
\end{table}

\begingroup
\section{Conclusions}
\label{sec:conclusions}

We have presented an automated shape optimisation framework for geometrically
nonlinear shell models on faceted triangulations. In particular, we considered
a five-parameter Naghdi model, although the same framework is applicable to
different models with only minor modifications. In contrast to isogeometric
approaches, our formulation does not require analytic or spline charts.
Instead, the geometry recovery and full-load state solve are recorded on the
automatic differentiation tape to evaluate first- and second-order shape
derivatives. This is a significant advantage, as it facilitates model
comparison and experimentation.

The framework is implemented in Firedrake and uses the shape optimisation
toolbox Fireshape and the optimisation library ROL to identify geometric
stiffening designs. Numerous numerical experiments and two industrial
benchmarks confirm the correctness and accuracy of the approach, in terms
of geometry recovery, nonlinear analysis, and design discovery. For example,
the optimised semi-cylinder design reduces the squared-displacement functional
by \semicylReductionPct\,\si{\percent}, and the optimised bracket design reduces
the elastic strain energy by \bracketReductionPct\,\si{\percent}. The bracket
develops the benchmark corrugation mechanism with a displacement allowance and
an achieved area increase comparable to those of the COMSOL benchmark, despite
using different shell, design-space and constraint formulations. On this test
case, we also computed optimal designs using linearised kinematics.
Interestingly, the framework successfully identified geometric stiffening
mechanisms despite significantly underestimating the displacement of the
initial geometry.

As expected, in light of the nonconvexity of the problem and the use of
gradient-based optimisation algorithms, discretisation parameters can
significantly affect the attainable designs. For instance, choosing a coarser
piecewise-linear control with a Helmholtz metric helps suppress
short-wavelength updates. Note that this still allows computations on an
independently refined state mesh. Including an area penalty reduces the
tendency to add closely spaced corrugations within a displacement bound.
Similarly, different initial designs can lead to different stiffening
mechanisms with varying performance.

The framework can be further extended both in terms of discretisation and
applications. On the one hand, it would be of interest to extend the geometry
recovery procedure to higher-order meshes, which would require revisiting its
mathematical foundations in addition to adapting its implementation. On the
other hand, the optimisation problems considered here do not account for
constraints on stress, stability, self-contact or manufacturability, which
would need to be properly formulated and incorporated into the optimisation
framework. Finally, extending the framework to elastoplastic shell models would
be of interest to investigate how material yielding and permanent deformation
influence the optimal stiffening mechanisms.

\par
\endgroup

\section*{CRediT authorship contribution statement}
\textbf{Ado Farsi:} Conceptualisation, Methodology, Software, Validation,
Investigation, visualisation, Writing (original draft).
\textbf{Alberto Paganini:} Writing (review and editing), Funding acquisition.

\section*{Declaration of competing interest}
The authors declare that they have no known competing financial interests or
personal relationships that could have appeared to influence the work reported in
this paper.

\section*{Code and data availability}
The supplementary reproducibility archive contains the example drivers,
exact input meshes, run arguments, selected numerical evidence, final controls,
figure inputs and checksum manifests. It includes explicit commands for the
reported baselines and the dependency versions recorded by the runs.
The implementation uses the open-source finite-element library Firedrake
(\url{https://github.com/firedrakeproject/firedrake}) together with the
shape-optimisation toolbox Fireshape
(\url{https://github.com/fireshape/fireshape}). The manuscript sources and
instructions for generating the clean manuscript, marked comparison and
reviewer response are included in the archive.

\section*{Acknowledgements}
Ado Farsi's work on this project was supported by the Impact Development Fund
(2025/26) of the University of Leicester.

\end{document}